\documentclass[12pt,english]{article}
\usepackage[margin=0.75in]{geometry}
\usepackage[utf8]{inputenc}
\usepackage[english]{babel}
\usepackage{amsthm}
\usepackage{graphicx}
\usepackage{url}
\usepackage{hyperref}
\usepackage{amsmath}
\usepackage{listings}
\usepackage{amsfonts}
\usepackage{xcolor}
\usepackage{textcomp}
\usepackage{caption}
\usepackage{subcaption}
\usepackage{xcolor}
\usepackage{ulem}

\begin{document}

\title{On the series expansion of the secondary zeta function}

\author{Artur Kawalec}

\date{}
\maketitle

\begin{abstract}
In article, we explore the secondary zeta function $Z(s)$, which is defined as a generalized zeta type of series over imaginary parts of non-trivial zeros of the Riemann zeta function $\zeta(s)$. This function has been analytically continued as a meromorphic function in $\mathbb{C}$ with one double pole and an infinity of simple poles.  The secondary zeta function is of interest because it can naturally represent an analytical formula for non-trivial zeros of the Riemann zeta function that we will explore, and we show that the non-trivial zeros can be generated directly from primes by introducing a new form of an explicit formula written in terms of the prime zeta function. Additionally, we will also give several new series expansions for $Z(s)$ and numerically compute these coefficients to high precision, and also develop several new methods to analytically extend $Z(s)$ to larger domains and develop algorithms to compute them.
\end{abstract}

\section{Introduction}
Let $t_n$ be the $n^{th}$ imaginary part of non-trivial zeros $\rho_n=\frac{1}{2}+it_n$ ($n=1,2,3...$) of the Riemann zeta function $\zeta(s)$ on the critical-line (assuming RH). The first few zeros have the values $t_1 = 14.13472514...$, $t_2 = 21.02203964...$, $t_3 = 25.01085758...$, and so on. Then, the secondary zeta function is defined as a generalized zeta series over these imaginary parts

\begin{equation}\label{eq:20}
Z(s)=\sum_{n=1}^{\infty}\frac{1}{t_n^{s}}
\end{equation}
for complex variable $s=\sigma+it$, as defined by Delsarte [10], Chakravarty [6][7][8]. The secondary zeta function (1) converges absolutely for $\Re(s)>1$ and can be analytically continued to the whole complex plane, where it is a meromorphic function with one double pole at $s=1$, and an infinity of simple poles at negative odd integers. This function has been very rarely studied in the literature, and recently appears in the works such as of Voros [23][24] and Arias De Reyna [2], who developed more theory to access this function, but the idea actually goes back to Mellin (translated by Voros [24, p.139]) who gave many kinds of series over the non-trivial roots, and even more so it goes back to Riemann who wrote down the first sum (1) for $s=2$ as pointed out by Arias De Reyna/R.P. Brent [3]. Although Chakravarty defined two secondary functions linked by a functional equation; the other function being related to certain sums over primes as in the Weil's explicit formula [25], or Guinand [12], which will not be discussed in this article, but the first kind (1) is the main focus of this study.

The main motivation for this article is that (1) provides a means of generating a true analytical formula for non-trivial zeros of the $\zeta$ function under a certain excitation as $s\to\infty$. In fact, we will give a new kind of explicit formula whereby each non-trivial zero can be individually computed from primes, thus essentially, transforming all primes into individual non-trivial zeros analytically.

We will further numerically explore the secondary zeta function for other values of $s$, and develop it into several different power series expansions, and also explore new methods to analytically continue it. We will perform extensive numerical computation of these expansion coefficients and verify these formulas numerically. We also remark that all of these results pose an interesting computational problem, as they are very difficult to compute to high precision.

The point of departure is a special case for the even values of $Z(2m)$ for $m\geq 1$ that we investigate next, where $m$ is a positive integer. If we consider the Hadamard product for the Riemann zeta function

\begin{equation}\label{eq:20}
\zeta(s)=\frac{\pi^{s/2}}{2(s-1)\Gamma(1+\frac{s}{2})}\prod_{\rho}^{}\left(1-\frac{s}{\rho}\right)
\end{equation}
and relating to the Riemann xi function gives the form

\begin{equation}\label{eq:20}
\xi(s) =\frac{(s-1)\Gamma(1+\frac{s}{2})}{\pi^{s/2}}\zeta(s)=\frac{1}{2}\prod_{\rho}^{}\left(1-\frac{s}{\rho}\right)
\end{equation}
satisfying the functional equation $\xi(s)=\xi(1-s)$. By taking the $m^{th}$ log-derivative of $\xi(\frac{1}{2}+it)$ at $t=0$ generates a formula for $Z(2m)$  as such

\begin{equation}\label{eq:20}
Z(2m)=-\frac{1}{2(2m-1)!}\frac{d^{(2m)}}{dt^{(2m)}}\log \xi(\frac{1}{2}+it)\Bigr\rvert_{t\to 0}
\end{equation}
for $m\geq 1$ and equating with the (lhs) of (3) we get the Voro's formula

\begin{equation}\label{eq:20}
\begin{aligned}
Z(2m) = (-1)^m \bigg[-\frac{1}{2(2m-1)!}\left[\log\zeta(\tfrac{1}{2})\right]^{(2m)}+\\
         -\frac{1}{4}\big((2^{2m}-1)\zeta(2m)+2^{2m}\beta(2m)\big)+2^{2m}\bigg]
\end{aligned}
\end{equation}
with the connecting identity

\begin{equation}\label{eq:20}
\begin{aligned}
\frac{(-1)^m}{2^m(m-1)!}\left[\log\Gamma(\tfrac{5}{4})\right]^{(m)} &=2^m\left[\frac{1}{2}\left(\left(1-2^{-m}\right)\zeta(m)+\beta(m)\right)-1\right]\\
\\
&=\frac{1}{2^m}\zeta(m,\tfrac{5}{4})=\sum_{n=1}^{\infty}\frac{1}{(\frac{1}{2}+2n)^m}
\end{aligned}
\end{equation}
as given by Voros [23][24]. The function $\zeta(s,a)$ is the classical Hurwitz zeta function and $\beta(s)$ is the Dirichlet beta function [1, p. 803], and the value for $\beta(2)$ is the Catalan's constant. The $\Gamma(s)$ function is the standard gamma function.  From (5), the first few $Z(2m)$ values can be generated as

\begin{equation}\label{eq:20}
\begin{aligned}
Z(2)&=\frac{1}{2}\left[\log\zeta(\tfrac{1}{2})\right]^{(2)}+\frac{1}{8}\pi^2+\beta(2)-4, \\
\\
Z(4)&=-\frac{1}{12}\left[\log\zeta(\tfrac{1}{2})\right]^{(4)}-\frac{1}{24}\pi^4-4\beta(4)+16,\\
\\
Z(6)&=\frac{1}{240}\left[\log\zeta(\tfrac{1}{2})\right]^{(6)}+\frac{1}{60}\pi^6+16\beta(6)-64, \\
\\
Z(8)&=-\frac{1}{10080}\left[\log\zeta(\tfrac{1}{2})\right]^{(8)}-\frac{17}{2520}\pi^8-64\beta(8)+256,\\
\\
Z(10)&=\frac{1}{725760}\left[\log\zeta(\tfrac{1}{2})\right]^{(10)}+\frac{31}{11340}\pi^{10}+256\beta(10)-1024,\\
\\
\end{aligned}
\end{equation}
and so on, where these values are described independently of the main non-trivial zero definition (1). The numerical values are computed in Table 3. This formula for the case $Z(2)$ was already known to Riemann (as pointed out by Arias De Reyna/R.P. Brent [3] as found in Riemann's Nacklass). Fortunately, these expressions are actually quite easy to compute to high accuracy, and will be used as a golden reference to calibrate other formulas for $Z(s)$ that we will study in this article.

When the sequence of zeros is ordered from smallest to largest such that

\begin{equation}\label{eq:20}
t_{1} < t_{2} < t_{3} < \ldots < t_{n},
\end{equation}
then the next consecutive reciprocal powers
\begin{equation}\label{eq:20}
O(t^{-s}_{n})\gg O(t^{-s}_{n+1})
\end{equation}
decay extremely fast as $s \to \infty$. Hence, one has the asymptotic leading term
 \begin{equation}\label{eq:20}
Z(s)=O\left(t_1^{-s}\right),
\end{equation}
which in turn can be used to extract the first zero in the limit as
\begin{equation}\label{eq:20}
t_{1} = \lim_{s\to\infty}\left[Z(s)\right]^{-\frac{1}{s}},
\end{equation}
and then, by successfully removing the previous zeros in the limit, the full recurrence formula for the $n^{th}+1$ zero becomes
\begin{equation}\label{eq:20}
t_{n+1} = \lim_{s\to\infty}\left(Z(s)-\sum_{k=1}^{n}\frac{1}{t_k^s}\right)^{-\frac{1}{s}},
\end{equation}
where they can be effectively removed from $Z(s)$ in the limit as $s\to\infty$. As a result, if we substitute (5) for $Z(s)$ into (12), then we obtain an analytical formula to generate non-trivial zeros
\begin{equation}\label{eq:20}
t_{n+1} = \lim_{m\to\infty}\left[\frac{(-1)^{m}}{2}\left(2^{2m}-\frac{1}{(2m-1)!}\left[\log\zeta(\tfrac{1}{2})\right]^{(2m)}-\frac{1}{2^{2m}}\zeta(2m,\tfrac{5}{4})\right)-\sum_{k=1}^{n}\frac{1}{t_{k}^{2m}}\right]^{-\frac{1}{2m}}
\end{equation}
as we have shown in [16][17][18], where we note that the main terms are only restricted to the even $2m$ derivatives of $\log \zeta(s)$ at $s=\frac{1}{2}$, and the Hurwitz zeta as a proper scaling factor (and $m$ is a positive integer).

To test this formula, we summarize numerical computations in Table $1$ for various limit values of $m$ from low to high, where we can observe the convergence to $t_1$ as $m$ increases. Already at $m=10$ we get several digits of $t_1$ as shown by last correct digit as underlined in blue color, and at $m=100$ we get well over $30$ digits. We performed even higher precision computations, and the result is clearly converging to $t_1$. A sample program written in Pari/Gp is given in Appendix A.

\begin{table}[hbt!]
\caption{The computation of $t_1$ by equation (13) for different limit variable $m$.} 
\centering 
\begin{tabular}{| c | c | c |} 
\hline 
\textbf{m} & $\boldsymbol{t_1}$ \textbf{(First 30 Digits)}  & \textbf{Significant Digits}\\  [0.5ex]
\hline 
$1$ & 6.578805783608427637281793074245 & 0  \\
\hline
$2$ & 12.806907343833847091925940068962 & 0 \\
\hline
$3$ & 13.809741306055624728153992726341 & 0 \\
\hline
$4$ & 14.038096225961619450676758199577 & 0 \\
\hline
$5$ & 14.\textcolor{blue}{\underline{1}}02624784431488524304946186056 & 1 \\
\hline
$6$ & 14.\textcolor{blue}{\underline{1}}23297656314161936112154413740 & 1 \\
\hline
$7$ & 14.1\textcolor{blue}{\underline{3}}0464459254236820197453483721 & 2 \\
\hline
$8$ & 14.1\textcolor{blue}{\underline{3}}3083993992268169646789606564 & 2 \\
\hline
$9$ & 14.13\textcolor{blue}{\underline{4}}077755601528384660110026302 & 3 \\
\hline
$10$ & 14.13\textcolor{blue}{\underline{4}}465134057435907124435534843 & 3 \\
\hline
$15$ & 14.1347\textcolor{blue}{\underline{2}}1950874675119831881762569 & 5 \\
\hline
$20$ & 14.13472\textcolor{blue}{\underline{5}}096741738055664458081219 & 6\\
\hline
$25$ & 14.13472514\textcolor{blue}{\underline{1}}055464326339414131271 & 9 \\
\hline
$50$ & 14.134725141734693\textcolor{blue}{\underline{7}}89641535771021 & 16 \\
\hline
$100$ & 14.134725141734693790457251983562 & 34 \\
\hline
\end{tabular}
\label{table:nonlin} 
\end{table}

We next demonstrate how to compute higher order non-trivial zeros. We set $m=250$ and compute
\begin{equation}\label{eq:20}
\begin{aligned}
Z(2m) = 7.18316934899718140841650578011166023417090863769600 & \\
     8517536818521464413577481501771580460474425539208\times 10^{-576}\ldots.
\end{aligned}
\end{equation}
to a very high precision. From this number, we extract the first $10$ non-trivial zeros, which are summarized in Table $2$ showing the last correct digit as underlined in blue. But in order to apply the recurrence formula, all the previous non-trivial zeros have to be known to an even higher precision than what was computed for $m$ (so for example we pre-computed previous zero to $1000$ decimal places in order to compute the $t_{n+1}$ zero) as one cannot use the same $t_n$ obtained earlier because it will a cause self-cancelation in (13), because the accuracy for $t_n$ must be much higher than $t_{n+1}$ in order to ensure convergence. So initially, for $m=250$ we started with an accuracy of $87$ digits after decimal place computed for $t_1$, and then it dropped to $7$ to $12$ digits by the time it gets to $t_{10}$ zero. There is also a sudden drop in accuracy when the gaps between non-trivial zeros get sporadically closer. Hence, these formulas are not very practical for computing high zeros as a super high numerical precision is required, for example, at the first Lehmer pair at $t_{6709}=7005.06288$, the gap between next zero is about $\sim 0.04$ would make this sepration even harder.  Also, the average gap between zeros gets smaller as $t_{n+1}-t_{n}\sim\frac{2\pi}{\log(n)}$, making the use of this formula progressively harder and harder to compute. But despite the numerical difficulty, these formulas constitute a true analytical representation of non-trivial zeros in closed-from, and without assuming any initial conditions.

\begin{table}[hbt!]
\caption{The $t_{n+1}$ computed by equation (13) for $m=250$.} 
\centering 
\begin{tabular}{| c | c | c | c|} 
\hline 
\textbf{n} & $\boldsymbol{t_{n+1}}$ & \textbf{(First 30 Digits)}  & \textbf{Significant Digits}\\  [0.5ex]
\hline 
$0$ & $t_{1}$ & 14.134725141734693790457251983562 & 87 \\
\hline
$1$ & $t_{2}$  & 21.022039638771554992628479593896 & 38 \\
\hline
$2$ & $t_{3}$ & 25.010857580145688763213790992562 & 43 \\
\hline
$3$ & $t_{4}$  & 30.424876125859513\textcolor{blue}{\underline{2}}09940851142395 & 16 \\
\hline
$4$ & $t_{5}$  & 32.9350615877391896906623689640\textcolor{blue}{\underline{7}}3 & 29 \\
\hline
$5$ & $t_{6}$  & 37.58617815882567125\textcolor{blue}{\underline{7}}190902153280 & 18 \\
\hline
$6$ & $t_{7}$  & 40.918719012147\textcolor{blue}{\underline{4}}63977678179889317 & 13 \\
\hline
$7$ & $t_{8}$  & 43.3270732809149995194961\textcolor{blue}{\underline{1}}7449701 & 22 \\
\hline
$8$ & $t_{9}$  & 48.005150\textcolor{blue}{\underline{8}}79831498066163921378664 & 7 \\
\hline
$9$ & $t_{10}$ & 49.77383247767\textcolor{blue}{\underline{2}}299146155484901550 & 12 \\
\hline
\end{tabular}
\label{table:nonlin} 
\end{table}

The key component of formula (13) is the

\begin{equation}\label{eq:1}
\left[\log\zeta(\tfrac{1}{2})\right]^{(2m)}
\end{equation}
even log-zeta derivative term, which can be computed to high precision by various numerical techniques, but they don't often give any more insights as to the nature of this term. And since the Hurwitz zeta in (13) doesn't contain any other information about non-trivial zeros, then this log-zeta even derivative is the true source of non-trivial zeros. So where do these zeros come from? We show actually that this term (15) can related to the prime numbers, and that shows how the apparent randomness of primes gets imprinted on the non-trivial zeros. In the previous article [19], we derived a general log-zeta formula

\begin{equation}\label{eq:1}
\left[\log\zeta(s)\right]^{(m)}=(-1)^m \frac{(m-1)!}{(s-1)^m}+\sum_{k=0}^{\infty}\frac{(s-a)^{k}}{k!}\left[[\log\zeta(a)]^{(k+m)}-\frac{(k+m-1)!}{(1-a)^{k+m}}\right]
\end{equation}
valid in some limited region for $Re(s)<1$ centered at $a$, where $a$ is a positive constant $1<a<c$, but is not too large as it is bounded by another constant $c=6.2$. And the radius of convergence is $R=a+2$. Hence by setting $s=\tfrac{1}{2}$ and $a=2$ and re-arranging we obtain the form

\begin{equation}\label{eq:1}
2^{2m}-\frac{1}{(2m-1)!}\big[\log\zeta(\frac{1}{2})\big]^{(2m)}= -\sum_{k=0}^{\infty}\frac{(-1)^k}{k!}(\frac{3}{2})^{k}\left[[\log\zeta(2)]^{(k+2m)}-(-1)^{k+2m}(k+2m-1)!\right]
\end{equation}
whereupon by inserting into (5) leads to an alternate expression for the secondary zeta function
\begin{equation}\label{eq:20}
\mathnormal{Z}(2m) =\frac{(-1)^{m}}{2}\left[\frac{1}{(2m-1)!}\sum_{k=0}^{\infty}\frac{(-1)^{k+1}}{k!}(\frac{3}{2})^k\Bigg(\big[\log\zeta(2)\big]^{(k+2m)}-(-1)^k(k+2m-1)!\Bigg)-\frac{1}{2^{2m}}\zeta(2m,\frac{5}{4})\right]
\end{equation}
as a function of the log-zeta derivatives at $a=2$ instead of at $a=\frac{1}{2}$, for which there is a whole lot of other interesting formulas to represent them, and one type in particular can be naturally expressed in terms of primes as such

\begin{equation}\label{eq:1}
\begin{aligned}
\left[\log\zeta(s)\right]^{(m)}&=(-1)^m\sum_{p}\sum_{n=1}^{\infty}\frac{1}{p^{ns}}n^{m-1}\log^{m}(p), \\
&=(-1)^m\sum_{n=2}^{\infty}\frac{\Lambda(n)}{\log(n)n^s}\log^{m}(n),\\
&=\sum_{n=1}^{\infty}P^{(m)}(ns)n^{m-1},\\
\end{aligned}
\end{equation}
all valid for $\Re(s)>1$, and where the von Mangoldt's function is defined as

\begin{equation}
\Lambda(n)= \left \{
\begin{aligned}
&\log p, &&\text{if}\ n=p^k \text{ for some prime and integer } k\geq 1 \\
&0 && \text{otherwise},
\end{aligned} \right.
\end{equation}
and $P(s)$ is the prime zeta function
\begin{equation}\label{eq:1}
P(s)=\sum_{p}\frac{1}{p^s}
\end{equation}
running over all primes $p=\{2,3,5,7\ldots\}$, and its derivatives are given by
\begin{equation}\label{eq:1}
P^{(m)}(s)=(-1)^m\sum_{p}\frac{\log^m(p)}{p^s}.
\end{equation}
As a result, by substituting the prime zeta as our favorite representation into (18) we obtain the explicit formula

\begin{equation}\label{eq:20}
\begin{aligned}
\mathnormal{Z}(2m) =\frac{(-1)^{m}}{2}\Bigg[\frac{1}{(2m-1)!}\sum_{k=0}^{\infty}\frac{(-1)^{k+1}}{k!}(\frac{3}{2})^k \Bigg(\sum_{j=1}^{\infty}P^{(k+2m)}(2j)j^{k+2m-1}&-(-1)^k(k+2m-1)!\Bigg)+\\
&-\frac{1}{2^{2m}}\zeta(2m,\frac{5}{4})\Bigg]
\end{aligned}
\end{equation}
for $m\geq 1$. Thus, the secondary zeta function (1) is expressed in term of the prime zeta derivatives for $2m$ even argument, and hence, the full recurrence formula for non-trivial zeros is manifest

\begin{equation}\label{eq:20}
\begin{aligned}
t_{n+1}=\lim_{m\to\infty}\Bigg[\frac{(-1)^{m}}{2}\Bigg(\frac{1}{(2m-1)!}\sum_{k=0}^{\infty}\frac{(-1)^{k+1}}{k!}(\frac{3}{2})^k \Bigg(\sum_{j=1}^{\infty}&P^{(k+2m)}(2j)j^{k+2m-1}-(-1)^k(k+2m-1)!\Bigg)+\\
&-\frac{1}{2^{2m}}\zeta(2m,\frac{5}{4})\Bigg)-\sum_{k=1}^{n}\frac{1}{t_{k}^{2m}}\Bigg]^{-\frac{1}{2m}}
\end{aligned}
\end{equation}
thus it describes a pure transformation of primes directly into non-trivial zeros, whereby each individual zero can be extracted in the limit. But we note that the number of primes in the prime zeta has to be extremely large (perhaps $\sim 10^{30}$) to achieve a few digits of accuracy for $t_1$ for higher $m$ limit value. When computing with these formulas numerically, we encourage the reader start with (13) first and then verify (18) using high precision numerical techniques, which should be easy, and finally transition to (24) using primes and carefully observe how fast it looses accuracy when using primes. Hence it is very impractical to compute the non-trivial zeros directly by (24) from primes, but however, high precision numerical techniques clearly verify that (24) is correct.

So far we've considered the special case for $Z(2m)$, where it naturally leads to a formula for non-trivial zeros, which has been the main motivation behind this article. But now we seek to explore $Z(s)$ for other values in $s$ domain, and so first we consider an analytical continuation as described by Delsarte [10], Chakravarty [6], and Arias De Reyna in [2] as:

\begin{equation}\label{eq:20}
Z(s)=\frac{1}{\Gamma(\frac{s}{2})}\int_{0}^{\infty}x^{\frac{s}{2}-1}\boldsymbol{\theta}(x) dx
\end{equation}
where $\Gamma(s)$ is the gamma function, and substituting a Jacobi type of series over the non-trivial roots

\begin{equation}\label{eq:20}
\begin{aligned}
\boldsymbol{\theta}(x)=\sum_{n=1}^{\infty}e^{-t^2_n x} = -\frac{1}{2\sqrt{\pi x}}\sum_{n=2}^{\infty}&\frac{\Lambda(n)}{\sqrt{n}}e^{-\frac{1}{4x}\log^2(x)}+e^{\frac{x}{4}}-\frac{\gamma+\log(16\pi^2 x)}{8\sqrt{\pi x}}+\\
&+\frac{1}{4\sqrt{\pi x}}\int_{0}^{\infty}e^{-\frac{u^2}{16x}}\left(\frac{1}{u}-\frac{e^{\frac{3}{4}u}}{e^u-1}\right)du
\end{aligned}
\end{equation}
which is still valid for $x>1$. By splitting the integral

\begin{equation}\label{eq:20}
\begin{aligned}
Z(s)=I_1(s)&+I_2(s)\\
\\
I_1(s)=\frac{1}{\Gamma(\frac{s}{2})}\int_{0}^{a}\boldsymbol{\theta}(x)x^{\frac{s}{2}-1}dx&,\quad I_2(s)=\frac{1}{\Gamma(\frac{s}{2})}\int_{a}^{\infty}\boldsymbol{\theta}(x)x^{\frac{s}{2}-1}dx
\end{aligned}
\end{equation}
for $a>0$, the lower integral $I_1(s)$ is singular at origin and is not entire, while the upper $I_2(s)$ is entire.  A rigorous analysis made by Arias de Reyna [2] analytically continues $I_1(s)$ and gives the following set of formulas by considering these parts separately:

\begin{equation}\label{eq:20}
I_1(s)=-P(s)+E(s)-S(s),
\end{equation}
where the prime term
\begin{equation}\label{eq:20}
P(s)=\frac{1}{2\sqrt{\pi}}\sum_{n=2}^{\infty}\frac{\Lambda(n)}{\sqrt{n}}\frac{\Gamma(\frac{1-s}{2},\frac{\log^2 n}{4a})}{\Gamma(\frac{s}{2})}\left(\frac{2}{\log n}\right)^{1-s}
\end{equation}
is entire, and the exponential term

\begin{equation}\label{eq:20}
E(s)=\frac{1}{\Gamma(\frac{s}{2})}\sum_{n=0}^{\infty}\frac{1}{4^n n!}\frac{a^{n+\frac{s}{2}}}{n+\frac{s}{2}}
\end{equation}
is also entire, and the singular term is

\begin{equation}\label{eq:20}
S(s)=a^{\frac{s-1}{2}}\frac{1}{\Gamma(\frac{s}{2})}\Bigg[\left(-\frac{2}{(s-1)^2}+\frac{\gamma+\log(16\pi^2a)}{s-1}\right)+\sum_{n=1}^{N}\frac{B_n(\frac{3}{4})}{n!}\frac{(4\sqrt{a})^n\Gamma(\frac{n}{2})}{s+n-1}\Bigg]+O(a^{\frac{N+\sigma}{2}})
\end{equation}
as $N\to\infty$, is meromorphic with a double pole at $s=1$ and simple poles at negative odd integers. Also, $\gamma=0.5772156649\ldots$ is the Euler's constant. The second integral in (27) (also labeled as $A(s)$) can be transformed into a series

\begin{equation}\label{eq:20}
I_2(s) = A(s)= \sum_{n=1}^{\infty}\frac{\Gamma(\frac{s}{2},a t_n^2)}{\Gamma(\frac{s}{2})}\frac{1}{t_n^s},
\end{equation}
where $\Gamma(s,a)$ is an incomplete gamma function. Thus, together the equations (25)-(32) constitute a full analytical continuation of (1) to $\mathbb{C}\backslash\{1,-1,-3,-5,\ldots\}$, and we will refer to it as the (ADR) algorithm [2]. But in order to apply these formulas numerically in a computer program, one needs to have a small database of prime numbers to compute the von Mangoldt function, and another database of non-trivial zeros to compute $I_2(s)$, and another database of Bernoulli polynomials $B_n(\frac{3}{4})$. When these parameters are available, then $Z(s)$ can be computed, and the accuracy of the output will depend on how well the accuracy of these parameters are leveraged together. Furthermore, these equations have already been implemented in a Python library $\textbf{mpmath}$ as a $\textbf{secondzeta}$ function, which are freely available in [15]. We also were able to successfully implement the (ADR) algorithm in Pari/Gp software package [22], where we can independently compute $Z(s)$, as well as its derivatives very accurately, which in this article, will serve as a golden reference to aide in further development of several different series representations of $Z(s)$. A sample program in Pari is given in Appendix B. Furthermore, the $\textbf{derivnum}$ function in Pari can compute the derivatives by numerical differentiation techniques, for example, the simplest way is by the limit definition

\begin{equation}\label{eq:20}
Z^{\prime}(s)=\lim_{h\to 0}\frac{Z(s+h)-Z(s)}{h},
\end{equation}
as well as the higher order derivatives by similar higher order difference methods.

The first few special values of $Z(s)$ are computed using the (ADR) algorithm and are summarized in Table 3 by properly setting the parameters $a$ and $N$. It turns out that (ADR) algorithm can compute $Z(s)$ to a very high degree of accuracy (100 digits) with little effort for many cases, but with more effort of improving numerical precision we computed successfully to over $3000$ digits and beyond. In Appendix C, we just give the first few odd values $Z(2n+1)$ to $1000$ decimal places for reference. In Fig. 1 we show a general plot of $Z(s)$ similarly as in [2], and in Fig. 2-6, we generate a few more plots of various zoomed-in cross sections of $Z(s)$. We note that near the negative poles the asymptotic growth is very sharp.

\begin{table}[hbt!]
\caption{Summary of $Z(s)$ computed for the first few values ($30$ digits) using (ADR) for $a=0.015$ and $N=100$.} 
\centering 
\begin{tabular}{| c | c |} 
\hline 
$\boldsymbol{s}$ & $\boldsymbol{Z(s)}$ \\ [0.5ex]
\hline 
$-6$ &  -0.134765625 \\
\hline
$-4$ &   0.0234375 \\
\hline
$-2$ &  -0.28125 \\
\hline
$-1.5$ & 0.543192013090835468509500468007 \\
\hline
$-0.5$ & 0.785321481872238131428414667627 \\
\hline
$0$ &    0.875 \\
\hline
$0.25$ & 1.046292354793822380550774149561 \\
\hline
$0.5$ &  1.549059995596196967158137380839 \\
\hline
$0.75$ & 4.003304705990492198451272426242 \\
\hline
$1.5$ &  0.247759676684890683692669878485 \\
\hline
$2$ &    0.023104993115418970788933810430 \\
\hline
$3$ &    0.000729548272709704215875518569 \\
\hline
$4$ &    0.000037172599285269686164866262 \\
\hline
$5$ &    0.000002231188699502103328640628 \\
\hline
$6$ &    0.000000144173931400973279695381 \\
\hline
$7$ &    0.000000009675344542702350408719 \\
\hline
$8$ &    0.000000000663031680252990869873 \\
\hline
$9$ &    0.000000000045991912392894862969 \\
\hline
$10$ &   0.00000000000321366415061660121 \\
\hline
$11$ &   0.00000000000022556506251559664\\
\hline
\end{tabular}
\label{table:nonlin} 
\end{table}

While the values for $Z(s)$ for $\Re(s)>1$ are of interest, there are other special values for $\Re(s)<1$, for example the half value has the representation

\begin{equation}\label{eq:20}
Z(\frac{1}{2})=-\frac{1}{2}\zeta(\frac{1}{2},\frac{5}{4})+\frac{1}{\pi\sqrt{2}}\int_{0}^{\infty}\frac{1}{\sqrt{t}}\left(\frac{\zeta^{\prime}}{\zeta}(\frac{1}{2}+t)+\frac{1}{t-\frac{1}{2}}\right)dt
\end{equation}
as can be generated from a formula in Voros [23, p. 683], which we further explore in Section 5. The value at origin is

\begin{equation}\label{eq:20}
Z(0)=\frac{7}{8},
\end{equation}
and values at negative even integers have a general closed-form

\begin{equation}\label{eq:20}
Z(-2m)=(-1)^{m}2^{-2m}(1-\frac{1}{8}E_{2m})
\end{equation}
with $E_{2m}$ being the Euler numbers as given in [7][23][24]. This results in rational values for negative even integer argument, for example, the first few values are generated as

\begin{equation}\label{eq:20}
\begin{aligned}
Z(-2)&=-\frac{9}{32},\quad\quad Z(-4)=\frac{3}{128},\quad\quad Z(-6)=-\frac{69}{512},\\
\\
Z(-8)&=-\frac{1377}{2048}, \quad\quad Z(-10)=-\frac{50529}{8192}, \quad\quad Z(-12)=-\frac{2702757}{32768},\\ \nonumber
\\
\end{aligned}
\end{equation}
and so on, also as shown in Table 3.

Henceforth in the following sections, we will further develop $Z(s)$ into several power series representations and provide very high-precision computation of these expansion coefficients.  Additionally, we also develop several new methods to analytically continue $Z(s)$ to a global domain and review the computational results.

\begin{figure}[!hbp]
  \centering
  \includegraphics[width=170mm]{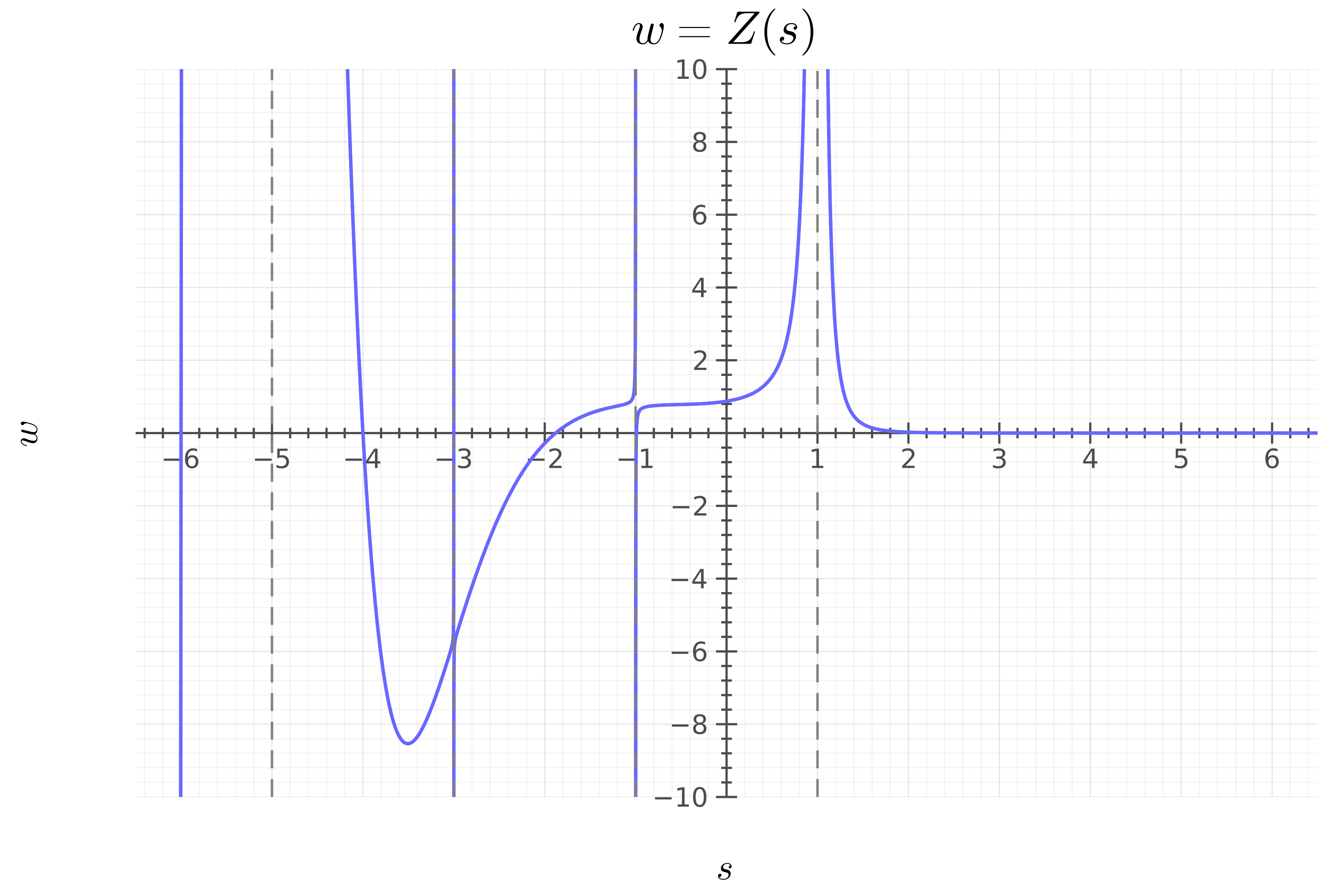}\\
  \caption{A plot of $w=Z(s)$ for $s\in (-6.5,6.5)$ }\label{1}
\end{figure}

\begin{figure}[!htbp]
  \centering
  \includegraphics[width=150mm]{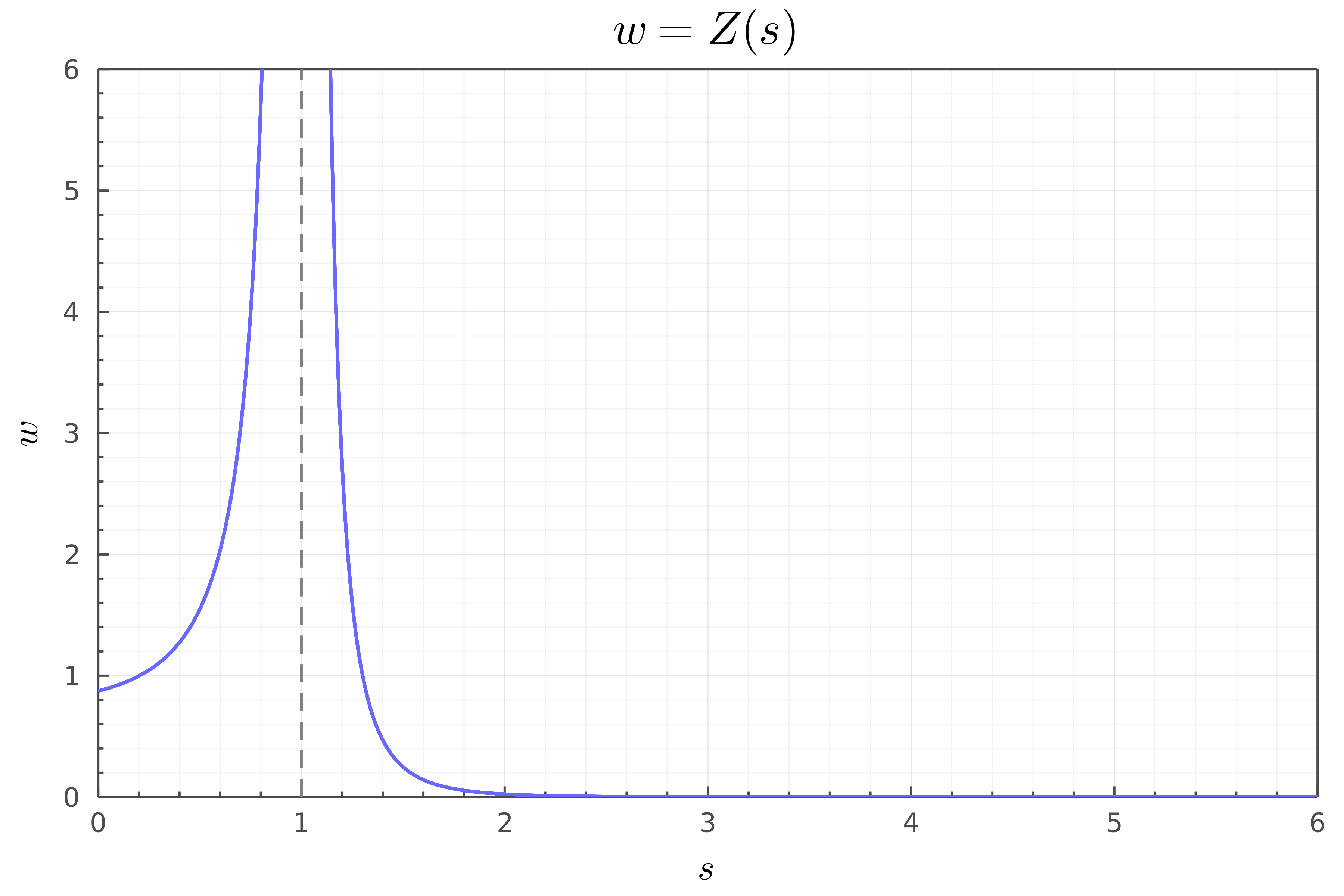}\\
  \caption{A zoomed-in plot of $w=Z(s)$ for $s\in (0,6)$ }\label{1}
\end{figure}

\begin{figure}[!htbp]
  \centering
  \includegraphics[width=150mm]{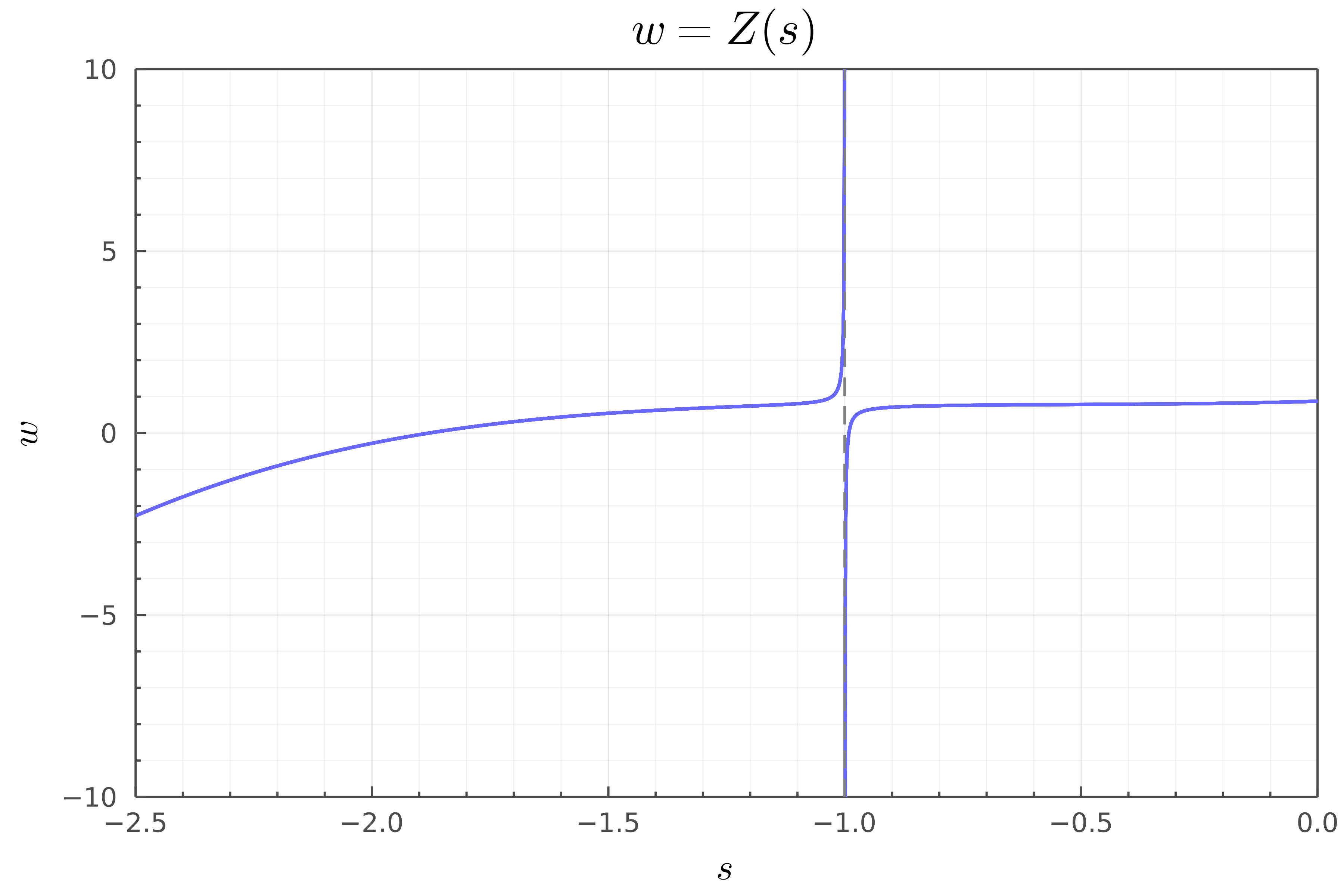}\\
  \caption{A zoomed-in plot of $w=Z(s)$ for $s\in (-2.5,0)$ }\label{1}
\end{figure}

\begin{figure}[!htbp]
  \centering
  \includegraphics[width=150mm]{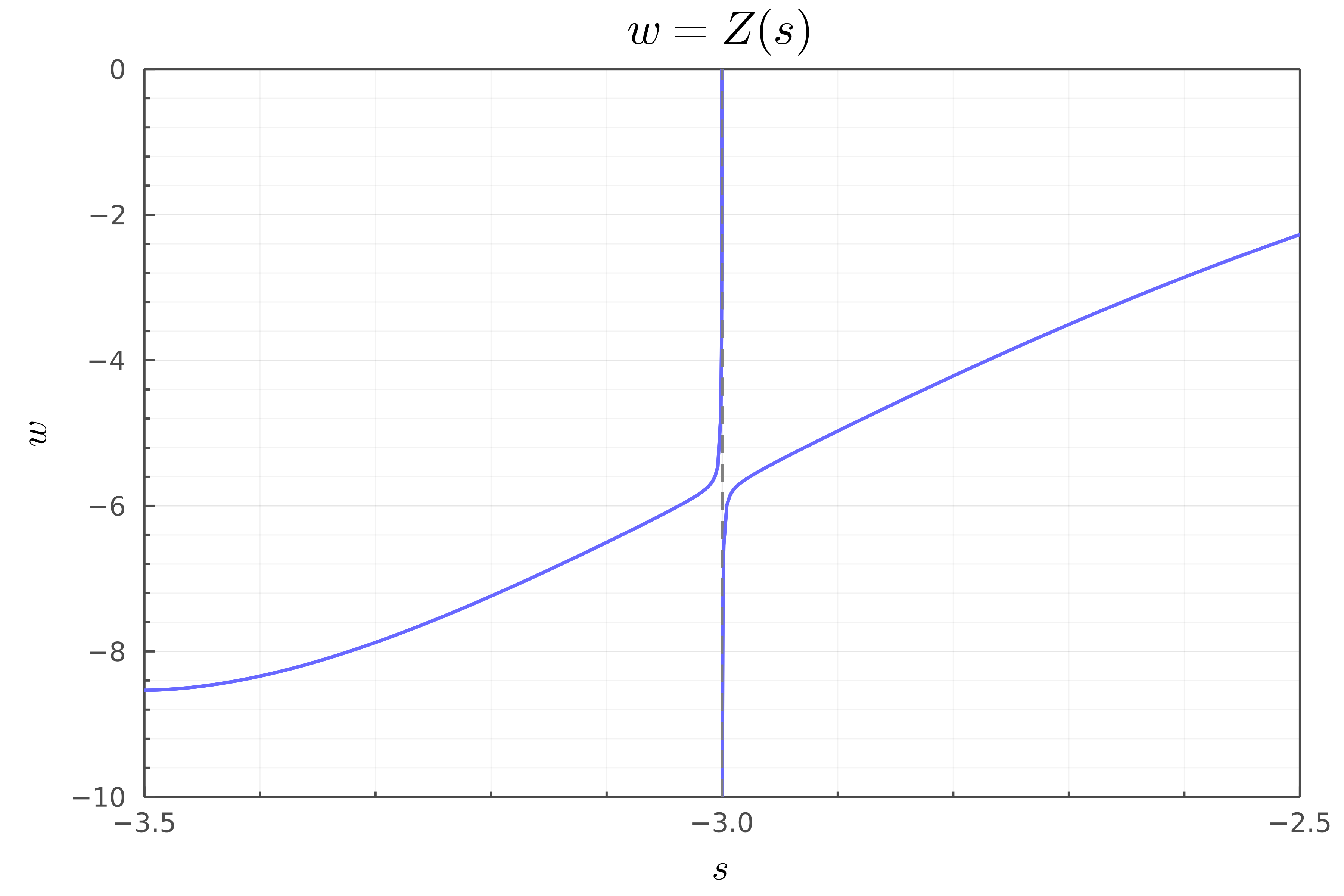}\\
  \caption{A zoomed-in plot of $w=Z(s)$ for $s\in (-3.5,-2.5)$ }\label{1}
\end{figure}

\begin{figure}[!htbp]
  \centering
  \includegraphics[width=150mm]{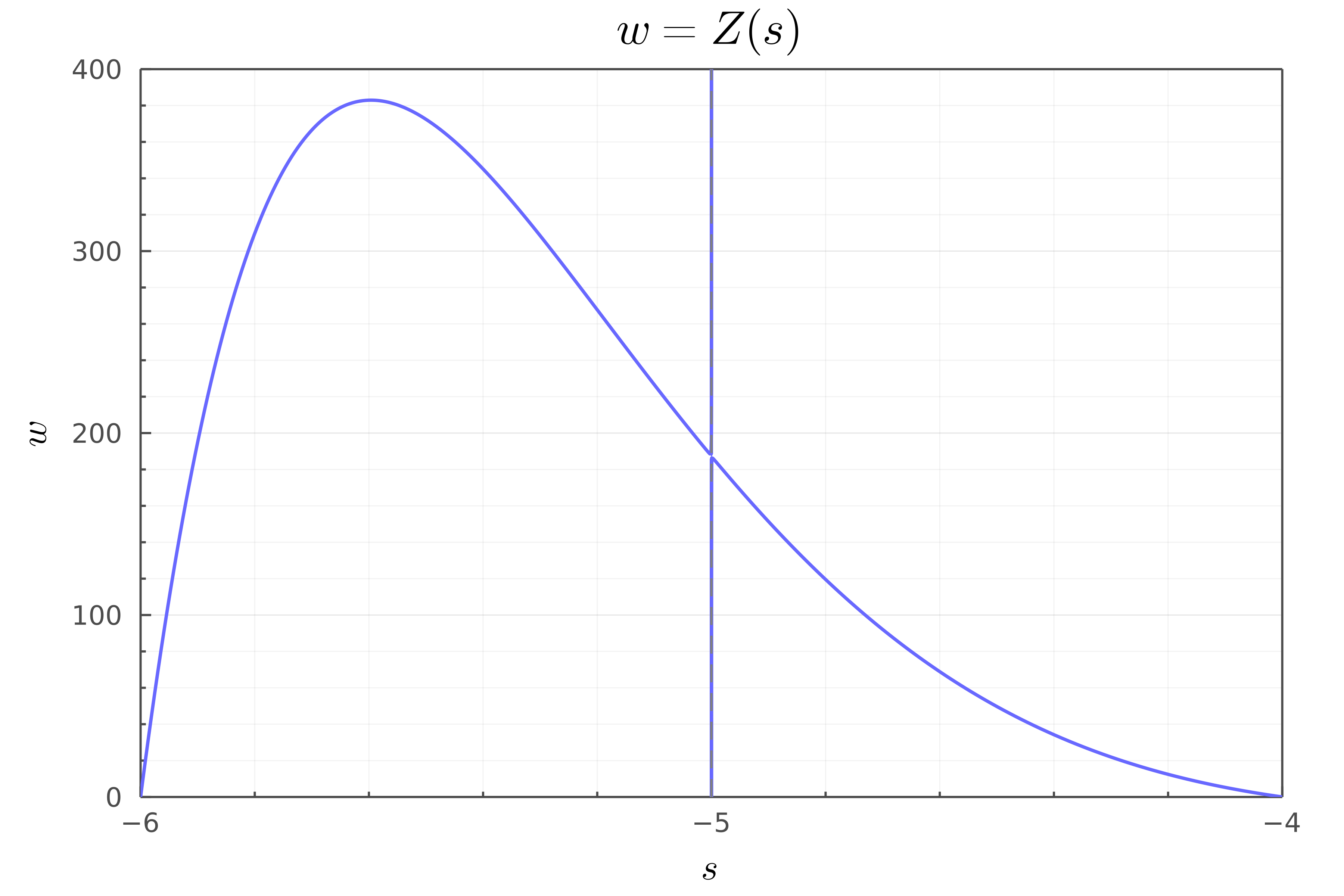}\\
  \caption{A zoomed-in plot of $w=Z(s)$ for $s\in (-6,-4)$ }\label{1}
\end{figure}

\begin{figure}[!htbp]
  \centering
  \includegraphics[width=150mm]{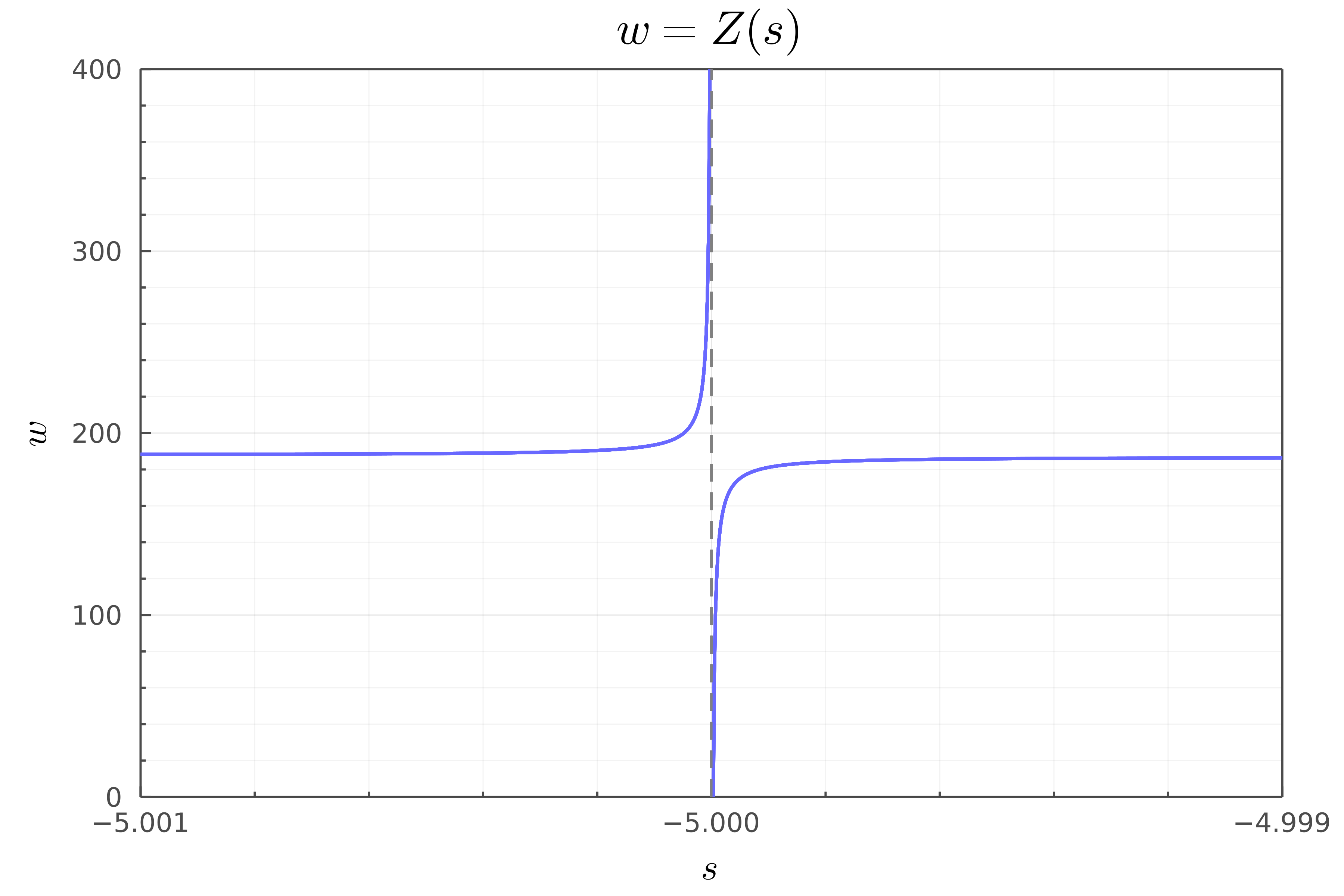}\\
  \caption{A zoomed-in plot of $w=Z(s)$ for $s\in (-5.001,-4.999)$ }\label{1}
\end{figure}

\newpage

\section{Simple Analytical Continuation of $Z(s)$}
In this section, we develop a very simple analytical continuation of (1) directly in terms of non-trivial zeros. Because (1) is defined for $\Re(s)>1$, we can write its Taylor series expansion about $a>1$ as

\begin{equation}\label{eq:20}
Z(s)=\sum_{m=0}^{\infty}Z^{(m)}(a)\frac{(s-a)^m}{m!},
\end{equation}
where the derivatives can be easily obtained from the original definition

\begin{equation}\label{eq:20}
Z^{(m)}(s)=(-1)^m\sum_{n=1}^{\infty}\frac{\log^m (t_n)}{t^s_n}
\end{equation}
for $\Re(s)>1$. However, the radius of convergence is only limited by distance from $a$ to the closest pole which exists at $s=1$, hence $R=a-1$. Now because there is a double pole at $s=1$, we will consider a new series expansion of the function instead

\begin{equation}\label{eq:20}
(s-1)^2 Z(s)=\sum_{m=0}^{\infty}A^{(m)}(a)\frac{(s-a)^m}{m!}
\end{equation}
that cancels the double pole, and $A_n(a)$ are the new generated expansion coefficients for $a>1$. Now, these $m^{th}$ derivatives can be expanded as

\begin{equation}\label{eq:20}
A^{(m)}(s)=\left[(s-1)^2Z(s)\right]^{(m)}
\end{equation}
which results in
\begin{equation}\label{eq:20}
A^{(m)}(a)=m(m-1)Z^{(m-2)}(a)+2m(a-1)Z^{(m-1)}(a)+(a-1)^2Z^{(m)}(a)
\end{equation}
for $m>1$, where it is implied that $A_0=(a-1)Z(a)$ for $m=0$ and $A_1=2(a-1)Z(a)+(a-1)^2Z^{(1)}(a)$ for $m=1$. By substituting the original definition (38) into (41), we obtain the formula

\begin{equation}\label{eq:20}
A^{(m)}(a)=(-1)^m\sum_{n=1}^{\infty}\frac{1}{t_n^a}\left[m(m-1)\log^{m-2}(t_n)-2m(a-1)\log^{m-1}(t_n)+(a-1)^2\log^m(t_n)\right]
\end{equation}
as a function of non-trivial zeros for $a>1$. As a result, by carrying over the double pole cancelation term to the (rhs) of (39), we obtain an analytical continuation as

\begin{equation}\label{eq:20}
Z(s)=\frac{1}{(s-1)^2 }\sum_{m=0}^{\infty}A^{(m)}(a)\frac{(s-a)^m}{m!}
\end{equation}
with an increased radius of convergence to $R=a+2$ due to being now limited by the next closest simple pole at $s=-1$. Hence, we can now access some small limited region below the line $\Re(s)<1$. So by substituting (42) into (37) we get the complete formula

\begin{equation}\label{eq:20}
\begin{aligned}
Z(s)=\frac{1}{(s-1)^2 }\sum_{m=0}^{\infty}&\sum_{n=1}^{\infty}\frac{(-1)^m}{t_n^a}\frac{(s-a)^m}{m!}\Big[m(m-1)\log^{m-2}(t_n)+\\
&-2m(a-1)\log^{m-1}(t_n)+(a-1)^2\log^m(t_n)\Big]
\end{aligned}
\end{equation}
which is now analytically continued to a larger domain centered at $a>1$ with $R=a+2$.
For example by setting $s=0$ and $a=2$, then we get a straightforward identity for $Z(0)$ value

\begin{equation}\label{eq:20}
Z(0)=\sum_{m=0}^{\infty}\sum_{n=1}^{\infty}\frac{2^m}{m!t_n^2}\Big[m(m-1)\log^{m-2}(t_n)-2m\log^{m-1}(t_n)+\log^m(t_n)\Big]=\frac{7}{8}
\end{equation}
where we obtain these explicitly from the non-trivial zeros.

\begin{table}[hbt!]
\caption{The computation of derivatives comparing (ADR) vs. Equation (38)} 
\centering 
\begin{tabular}{| c | c | c |} 
\hline
\textbf{m} & \textbf{$\boldsymbol{Z^{(m)}(2)}$ (ADR) Algorithm} & \textbf{$\boldsymbol{Z^{(m)}(2)}$ $\boldsymbol{10^6}$ non-trivial zeros equation (38)}\\
\hline 
$0$ &   0.02310499311541897078893381043 & 0.0231\textcolor{blue}{\underline{{0}}}1687559638772401176500051  \\
\hline
$1$ &  -0.09262185134364910430268093267 &-0.09\textcolor{blue}{\underline{{2}}}574299784035646356852860502 \\
\hline
$2$ &   0.41435919734495161051196316435 & 0.4\textcolor{blue}{\underline{{1}}}3671337197827754632029431426 \\
\hline
$3$ &  -2.09748326735062596330899544392 &-2.\textcolor{blue}{\underline{{0}}}87470336492663460380267586032 \\
\hline
$4$ &   12.0980012213366698631646223341 & 11.95120633944930957890920147222 \\
\hline
$5$ &  -79.4863416449377472053529935041 &-77.31676991503012640740477769007 \\
\hline
$6$ &   591.484015193663096366355496015 & 559.1210836688477955871718090022 \\
\hline
$7$ &  -4942.74891352383061750707285662 &-4454.867705330352429681931379267 \\
\hline
$8$ &   45959.9523623231656529494440178 & 38515.18899427343013017759525815 \\
\hline
$9$ &  -471395.877139293722056777184099 &-356197.9730212180865627942153289 \\
\hline
$10$ &  5291504.90926829401845873416044 & 3480169.820209129396024471479426 \\
\hline
\end{tabular}
\label{table:nonlin} 
\end{table}

But verifying these formulas numerically is much harder as it takes a very large database of non-trivial zeros to even get several digits. To demonstrate this, in Table 4 we compute the $n^{th}$ derivatives of $Z(s)$ about $a=2$ using the (ADR) algorithm and also by the original derivative formula (38) using a database of $1$ million non-trivial zeros. We observe that the accuracy of formula (38) starts to degrade very quickly as shown by the last correct underlined blue decimal place. As it is seen, for $m=3$ we're quickly down to $1$ decimal place. Because of that, our database of non-trivial zeros is clearly insufficient to compute $Z(s)$ to high enough accuracy, so instead, we will solely rely on the (ADR) algorithm, where we can compute $Z^{(m)}(a)$ to a super high precision independently of the non-trivial zeros. In Table 5, we compute the $A_m(a)$ coefficients up to $50$ decimal places to high precision, and using this database, we verify equation (43) for $Z(s)$ for several different values.

\begin{table}[hbt!]
\caption{The first $50$ $A^{(n)}(a)$ expansion coefficients generated about $a=2$ using the (ADR) algorithm.} 
\centering 
\begin{tabular}{| c | c |} 
\hline
$\boldsymbol{n}$  & $\boldsymbol{A^{(n)}(2)}$ \textbf{(also in 50 digits)}\\ [0.5ex]
\hline
0  & 0.0231049931154189707889338104303390140033817603974 \\
\hline
1  &-0.0464118651128111627248133118129024105529088846306 \\
\hline
2  & 0.0900817782011931348791070545201424412035427300596 \\
\hline
3  &-0.1670591913428109260533020538390815940736732110748 \\
\hline
4  & 0.2904454506710814828362167550613970584089615123387 \\
\hline
5  &-0.4559947785835678398866790408212000947406380442043 \\
\hline
6  & 0.5879520945102257970582439904890879126003237537305 \\
\hline
7  &-0.3990498999326510029216395829754796695201463433312 \\
\hline
8  &-0.9254032129908306478139112935193407118538981348969 \\
\hline
9  &5.34360880745523580356254390100629478500194977023364 \\
\hline
10 &-16.920908495513911359559945349109370729705587993182 \\
\hline
\vdots & \vdots \\
\hline
20 &-6.0687262197471062787446495667638723078912968550157$\times 10^{6}$ \\
\hline
30 &-1.1391213955114704007153966127464503426710621074724$\times 10^{16}$  \\
\hline
40 &-5.9339311931520365805736015140877173066785041127093$\times 10^{26}$  \\
\hline
50 &-3.7459309104238050373812710645367898060134580322909$\times 10^{38}$  \\
\hline
\end{tabular}
\label{table:nonlin} 
\end{table}

To verify these coefficients, we compute $Z(s)$ for several known test values. And as before, the last correct digit is underlined from the raw output. The value at origin is computed as
\begin{equation}\label{eq:20}
\begin{aligned}
Z(0)&\approx \sum_{n=0}^{50}A_n(2)\frac{(-2)^n}{n!}\\
&\approx0.875000000\underline{0}27734137822908541582\ldots,
\end{aligned}
\end{equation}
which is accurate to $10$ decimal places thus verifying (45), and another is computed as
\begin{equation}\label{eq:20}
\begin{aligned}
Z(\frac{1}{2})&\approx4\sum_{n=0}^{50}A_n(2)\frac{(-\frac{3}{2})^n}{n!}\\
&\approx 1.549059995596196\underline{9}98571023123734\ldots,
\end{aligned}
\end{equation}
which is accurate to $16$ decimal places, and another is computed as
\begin{equation}\label{eq:20}
\begin{aligned}
Z(3)&\approx \frac{1}{4}\sum_{n=0}^{50}A_n(2)\frac{1}{n!}\\
&\approx0.0007295482727097042158755\underline{1}7799\ldots,
\end{aligned}
\end{equation}
which is accurate to $26$ decimal places, where we observe the correct convergence and to a reasonable accuracy. In this expansion, this series is centered about $a=2$ having a radius of convergence $R=3$, so we cannot compute outside of that circle, but the importance is that we can already work in some limited region $\Re(s)<1$ with an only knowledge of $Z(s)$ for $\Re(s)>1$.

In order to further analytically continue $Z(s)$, we have remove more singularities, so we next will use the gamma function to cancel all poles at negative odd integers. Hence we consider the new function instead

\begin{equation}\label{eq:20}
\begin{aligned}
B(s)&=\frac{1}{\Gamma\Big(\dfrac{s+1}{2}\Big)}A(s) \\
\\
    &=\frac{1}{\Gamma\Big(\dfrac{s+1}{2}\Big)}(s-1)^2 Z(s)
\end{aligned}
\end{equation}
which is entire. Then, similarly as before, we write its Taylor series expansion

\begin{equation}\label{eq:20}
B(s)=\sum_{n=0}^{\infty}B^{(n)}(a)\frac{(s-a)^n}{n!}
\end{equation}
about $a>1$ again. The $n^{th}$ derivatives are generated as

\begin{equation}\label{eq:20}
B^{(n)}(s)=\Bigg[\frac{1}{\Gamma\Big(\dfrac{s+1}{2}\Big)}A(s)\Bigg]^{(n)}.
\end{equation}
The first two simple cases are
\begin{equation}\label{eq:20}
B^{}(s)=\frac{1}{\Gamma\Big(\dfrac{s+1}{2}\Big)}A(s),
\end{equation}
and the next one is
\begin{equation}\label{eq:20}
B^{\prime}(s)=\frac{1}{\Gamma\Big(\dfrac{s+1}{2}\Big)}A^{\prime}(s)-\frac{\Gamma^{\prime}\Big(\dfrac{s+1}{2}\Big)}{2\Gamma\Big(\dfrac{s+1}{2}\Big)^2}A(s)
\end{equation}
and so on. From (51) these $n^{th}$ derivatives expand to generate long such sequences due the chain rule, but numerically, they can be evaluated without a problem. We stress that the $A^{(n)}(a)$ can be computed from (1) directly, the main definition of the secondary zeta function when $a>1$, although it would be very difficult since a very large dataset of non-trivial zeros is required (perhaps one trillion) to get some reasonable numerical accuracy.

Now finally, we carry over the pole cancelation terms to the (lhs) of equation (49), we obtain a full analytical continuation

\begin{equation}\label{eq:20}
Z(s)=\frac{1}{(s-1)^2}\Gamma\Big(\dfrac{s+1}{2}\Big)\sum_{n=0}^{\infty}B^{(n)}(a)\frac{(s-a)^n}{n!}
\end{equation}
to $\mathbb{C}\backslash\{1,-1,-3,-5,\ldots\}$. This series can be centered at any $a$, but again, if $a>1$, then we can compute these expansion coefficients directly from the non-trivial zeros by equation (38), the main definition. But if $a<1$, then we will need different methods to compute these coefficients, and one method is of course the (ADR) algorithm, and in later sections, we will explore yet another independent method to generate these coefficients and without involving the non-trivial zeros.

\begin{table}[hbt!]
\caption{The first $50$ $B^{(n)}(a)$ expansion coefficients generated about $a=2$ by the (ADR) algorithm.} 
\centering 
\begin{tabular}{| c | c |} 
\hline
$\boldsymbol{n}$  & $\boldsymbol{B^{(n)}(2)}$ \textbf{(also in 50 digits)}\\ [0.5ex]
\hline
0 & 0.02607119288732401047712116299152947391229912695003 \\
\hline
1 &-0.05284585027436758872350351280594610669109818059365 \\
\hline
2 & 0.09747321487109608112019530365118296442929494602501 \\
\hline
3 &-0.15437090275753191135675664619278824153789839328444 \\
\hline
4 & 0.17654980634111359997175931414678567312778051279362 \\
\hline
5 &-0.00761229224736873022796626365155488797869442759740 \\
\hline
6 &-0.71721810717902008390828161105921240809381964913316 \\
\hline
7 & 2.55329742785210396006527280641471520620688162418726 \\
\hline
8 &-5.48123464868420465662992761530694912102945680986895 \\
\hline
9 & 5.59465678718268556542139228571337332842747630566711 \\
\hline
10& 14.1947998112068264846977684784000007495275789004403 \\
\hline
\vdots & \vdots \\
\hline
20& 2.16833981007094601718623725913192832212063088911808$\times 10^{6}$ \\
\hline
30& 7.70865253176087393689549717026595284167063191529794$\times 10^{11}$ \\
\hline
40& 3.73226885559764450638821041985553610058096587866085$\times 10^{17}$ \\
\hline
50&-4.92014848966196918922153672240236366804055822236138$\times 10^{23}$ \\
\hline
\end{tabular}
\label{table:nonlin} 
\end{table}

Next, in Table 6, we compute the first $B^{(n)}(a)$ derivatives centered at $a=2$ using equation (51), where we've reused all the 50 values of $A^{(n)}(a)$ from Table 5.  And then, we compute $Z(s)$ based on the equation (54) for the first few values. We compute

\begin{equation}\label{eq:20}
\begin{aligned}
Z(-4)&\approx\frac{1}{25}\Gamma(-\frac{3}{2})\sum_{n=0}^{50}B_n(2)\frac{(-6)^n}{n!}\\
&\approx 0.0\underline{2}0885048364271929740437763122\ldots,
\end{aligned}
\end{equation}
which is accurate to $2$ decimal places. The fewer accurate decimals means that we need just more terms of the Taylor series to get better accuracy in this region,  we next compute
\begin{equation}\label{eq:20}
\begin{aligned}
Z(-2)&\approx\frac{1}{9}\Gamma(-\frac{1}{2})\sum_{n=0}^{50}B_n(2)\frac{(-4)^n}{n!}\\
&\approx  -0.2812499999\underline{9}2678995349687792342\ldots,
\end{aligned}
\end{equation}
which is accurate to $11$ decimal places, and another is computed as
\begin{equation}\label{eq:20}
\begin{aligned}
Z(0)&\approx\sqrt{\pi}\sum_{n=0}^{50}B_n(2)\frac{(-2)^n}{n!}\\
&\approx 0.8749999999999999999999999\underline{9}0510\ldots,
\end{aligned}
\end{equation}
which is accurate to $26$ decimal places, and another is computed as
\begin{equation}\label{eq:20}
\begin{aligned}
Z(3)&\approx\frac{1}{4}\sum_{n=0}^{50}B_n(2)\frac{1}{n!}\\
&\approx 0.0007295482727097042158755185690939705033518\underline{1}596658\ldots,
\end{aligned}
\end{equation}
which is accurate to $44$ decimal places, and one more is computed as
\begin{equation}\label{eq:20}
\begin{aligned}
Z(5)&\approx\frac{1}{8}\sum_{n=0}^{50}B_n(2)\frac{3^n}{n!}\\
&\approx 0.00000223118869950\underline{2}319193060806\ldots,
\end{aligned}
\end{equation}
which is accurate to $18$ decimal places, but now the accuracy starts to drop again as we get farther away from the center $a$.  We quickly observe that this series works well in the limited region below the line $\Re(s)<1$. We also observe that the value for $Z(0)$ in equation (57) achieved $26$ decimal places, while the same value in equation (46) based on the $A^{(n)}(a)$ coefficients was accurate only to $10$ decimal places. No new computation was  made, the $B^{(n)}(a)$ coefficients were derived directly from $A^{(n)}(a)$ in Table $5$ by interweaving the derivatives with the $\Gamma$ function based on equation (51). This means that the removal of the simple poles at negative odd integers caused an improved convergence of the series, especially near those poles; hence the $B$ series (54) is actually preferred instead of $A$ series (43) because of the removal of all these singularities.

\section{Taylor series at origin}
In the previous section, we obtained a full Taylor series expansion about $a>1$ since it can be easily obtained from the original definition (1), and now, we'll investigate a similar series expansion about the origin when $a=0$ and write

\begin{equation}\label{eq:20}
Z(s)=\sum_{n=0}^{\infty}Z^{(n)}(0)\frac{s^n}{n!},
\end{equation}
which has radius of convergence $R=1$ due to again the double pole at $s=1$ and a second simple pole at $s=-1$. We already know that the value $Z(0)=\frac{7}{8}$, and the importance of this series is that there is a closed-form for the first derivative as

\begin{equation}\label{eq:20}
Z^{(1)}(0)=\frac{1}{2}\log\left(\frac{2^{\frac{11}{4}}\sqrt{\pi}}{\Gamma(\frac{1}{4})|\zeta(\frac{1}{2})|}\right)
\end{equation}
as shown by Voros [24, p. 85] (but we've added a factor of $\frac{1}{2}$ since Voros originally defined a re-scaled function $\mathcal{Z}(s)=Z(2s)$). We've also used this value extensively to calibrate the accuracy of the (ADR) algorithm.

And to get a full analytical continuation as we did in (53) we have

\begin{equation}\label{eq:20}
B^{(n)}(s)=\Bigg[(s-1)^2\frac{1}{\Gamma\Big(\dfrac{s+1}{2}\Big)}Z(s)\Bigg]^{(n)}.
\end{equation}
and at $s=0$ we get the full form
\begin{equation}\label{eq:20}
Z(s)=\frac{1}{(s-1)^2}\Gamma(\dfrac{s+1}{2})\sum_{n=0}^{\infty}B^{(n)}(0)\frac{s^n}{n!},
\end{equation}
which is valid in $\mathbb{C}\backslash\{1,-1,-3,-5,\ldots\}$. In Tables 7-8 we tabulate a high-precision computation of these coefficients for further study. The first two $B$ coefficients are therefore

\begin{equation}\label{eq:20}
B^{(0)}(0) =\frac{7}{8\sqrt{\pi}} \nonumber
\end{equation}
and from (53) we have

\begin{equation}\label{eq:20}
B^{(1)}(0) = -\frac{7}{4\sqrt{\pi}}-\frac{7}{16\pi}\Gamma^{\prime}(\frac{1}{2})+\frac{1}{2\sqrt{\pi}}\log\left(\frac{2^{\frac{11}{4}}\sqrt{\pi}}{\Gamma(\frac{1}{4})|\zeta(\frac{1}{2})|}\right) \nonumber
\end{equation}
and so on.

\begin{table}[hbt!]
\caption{The computation of $Z^{(n)}(0)$ (first 50 digits) by the (ADR) algorithm.} 
\centering 
\begin{tabular}{| c | c |} 
\hline
$\boldsymbol{n}$  & $\boldsymbol{Z^{(n)}(0)}$ \textbf{(also in 50 digits)}\\ [0.5ex]
\hline
0 & 0.875 \\
\hline
1 & 0.40590897213384756573995541397622554173832123404575 \\
\hline
2 & 1.56016236721756204469457369276272045935954269246077 \\
\hline
3 & 5.85079911124888264238220272801114910187113684144278 \\
\hline
4 & 25.2820311559571370258436434355729365053785809270153 \\
\hline
5 & 151.816823257781283090093014014169330878144109484301 \\
\hline
6 & 1005.91850911935465340848532909848677515234861530874 \\
\hline
7 & 7926.95848151906906780109595943735168033073981192661 \\
\hline
8 & 69281.0198568462452979234407900206825602825941344497 \\
\hline
9 & 686082.342172282824297860625726151573120546946894154 \\
\hline
10 & 7390380.6814175697837386165354747845013738393673307 \\
\hline
\vdots & \vdots \\
\hline
20 & 8.8268837140573609879990792592995333763054815128704$\times 10^{23}$ \\
\hline
30 & 1.3845347819988183483182022346624523558527997567154$\times 10^{33}$ \\
\hline
40 & 5.5573853743169003196882177991678993777611910516202$\times 10^{48}$ \\
\hline
50 & 2.5556286903717565654096115139313266385753855580946$\times 10^{65}$ \\
\hline
\end{tabular}
\label{table:nonlin} 
\end{table}

\begin{table}[hbt!]
\caption{The computation of $B^{(n)}(0)$ (first 50 digits) by the (ADR) algorithm.} 
\centering 
\begin{tabular}{| c | c |} 
\hline
$\boldsymbol{n}$  & $\boldsymbol{B^{(n)}(0)}$ \textbf{(also in 50 digits)}\\ [0.5ex]
\hline
0 & 0.49366588560428675107956952011567601261354430066287 \\
\hline
1 & -0.2736631993176991598064727716007833019766028144175 \\
\hline
2 & -0.6706720001793549880150016236038318899020204779574 \\
\hline
3 & 2.52188950096038926543249485356486551201807467724506 \\
\hline
4 & -3.4690248685075380539828994198565947646279961286361 \\
\hline
5 & -4.5690329917932599989769338540021897425280877471131 \\
\hline
6 & 40.6224037422036274157108952153073803251327066620718 \\
\hline
7 & -120.15747006396059169864205595849235341922888406894 \\
\hline
8 & 149.021561370514060320742352101688964090619215586130 \\
\hline
9 & 435.214123267853958975528589245925372198387078933914 \\
\hline
10 & -3347.828851792287018494762500504433926685952633985 \\
\hline
\vdots & \vdots \\
\hline
20 & -8.660364941727636469575190868981760507255293417496$\times 10^{8}$  \\
\hline
30 & -6.494091875636945962330127599787046043900599551078$\times 10^{14}$ \\
\hline
40 & -1.197691890793835431189886336725116607513908699743$\times 10^{21}$ \\
\hline
50 & -3.330875317972086567188644428857137654539498341328$\times 10^{27}$ \\
\hline
\end{tabular}
\label{table:nonlin} 
\end{table}

\newpage
\section{Laurent series at the double pole}
The principal part of the Laurent expansion of $Z(s)$ around the double pole is originally given by Delsarte [10] and Chakravarty [7] as

\begin{equation}\label{eq:20}
Z(s)=\frac{1}{2\pi(s-1)^2}-\frac{\log(2\pi)}{2\pi(s-1)}+O(1)
\end{equation}
and the regular part later by \text{Bondarenko-Ivi$\acute{c}$-Saksman-Seip} in [14] and also in R.P. Brent [5], to complete the form

\begin{equation}\label{eq:20}
Z(s)=\frac{1}{2\pi(s-1)^2}-\frac{\log(2\pi)}{2\pi(s-1)}+C_0+\sum_{n=1}^{\infty}C_n\frac{(s-1)^n}{n!}.
\end{equation}
This series is centered at $a=1$ with radius of convergence being $R=2$. The value of the constant $C_0$ has a special significance as being an analogue of the Euler's harmonic sum

\begin{equation}\label{eq:20}
\gamma=\lim_{k\to\infty}\Bigg\{\sum_{n=1}^{k}\frac{1}{n}-\log(k)\Bigg\}
\end{equation}
for the Euler's constant, where we have a similar harmonic sum over non-trivial roots

\begin{equation}\label{eq:20}
H =\lim_{k\to\infty}\Bigg\{\sum_{n=1}^{k}\frac{1}{t_n} - \frac{1}{4\pi}\log^2\left(\frac{t_k}{2\pi}\right)\Bigg\}
\end{equation}
where the limit is with respect to index $k$. This constant has been shown to exist and computed to $18$ decimal places by R.P Brent in [3][4] (without assuming RH) using a database of over $10^{10}$ non-trivial for the Brent's latest convergence acceleration method for $Z(s)$. This constant is related to the $0^{th}$ order expansion coefficient in (65) by

\begin{equation}\label{eq:20}
C_0=H+\frac{\log^2(2\pi)}{4\pi}
\end{equation}
and has also been studied and estimated by Hassani [13]. The reason for these two parts is that $H$ comes from the regular part of $Z(s)$ while $\frac{\log^2(2\pi)}{4\pi}$ comes from the principal part, and so they are merged together. From these results, we extract it in the limit as
\begin{equation}\label{eq:20}
C_0=\lim_{s\to 1}\Bigg\{Z(s)-\frac{1}{2\pi(s-1)^2}+\frac{\log(2\pi)}{2\pi(s-1)}\Bigg\}.
\end{equation}
Now using the machinery of the (ADR) that we saw in the previous sections, we compute

\begin{equation}\label{eq:20}
\begin{aligned}
H =-&0.01715940430709814945419161427396429261536011174\\
      &8183591215159125137067824458151406779819175728159\\
      &853433378444762\ldots\\
\end{aligned}
\end{equation}
to $111$ decimal places, which we believe to be very accurate assuming (RH). We note that even with the (ADR) algorithm we found that it was much harder to compute it near the double pole at $s=1$ than for other values. And in Table 9, we tabulate $C_0$ and the next higher order constants $C_n$ to high precision computed by

 \begin{equation}\label{eq:20}
C_n=\lim_{s\to 1}\frac{d^n}{ds^n}\Bigg[Z(s)-\frac{1}{2\pi(s-1)^2}+\frac{\log(2\pi)}{2\pi(s-1)}\Bigg].
\end{equation}

To verify these coefficients, we compute the test value

\begin{equation}\label{eq:20}
\begin{aligned}
Z(0)&\approx \frac{1}{2\pi}+\frac{\log(2\pi)}{2\pi}+\sum_{n=0}^{50}C_n\frac{(-1)^n}{n!}\\
&\approx 0.8750000000000000\underline{0}294495\ldots,
\end{aligned}
\end{equation}
using the coefficients from Table 9, and the result is accurate to $17$ decimal places. And similarly, the other special values near $1$ are also that accurate. This is a good confirmation of the validity of these coefficients.

\begin{table}[hbt!]
\caption{The computation of $C_n$ (first 50 digits) by the (ADR) algorithm.} 
\centering 
\begin{tabular}{| c | c |} 
\hline
$\boldsymbol{n}$  & $\boldsymbol{C_n}$ \textbf{(also in 50 digits)}\\ [0.5ex]
\hline
0 &  0.2516367513127059665334663293426453755147595873836 \\
\hline
1 & -0.1300444859118885707285274533988846777460553964263 \\
\hline
2 &  0.0824214912550528039526632284933172430791521350021 \\
\hline
3 & -0.0321581827282544905964296099391141952179545405019 \\
\hline
4 & -0.0531801364893419772868761573698112582469915802523  \\
\hline
5 &  0.2110321083617385257637243839874627961215847994456  \\
\hline
6 & -0.4933371057135871285817870279321636575675112589435  \\
\hline
7 &  0.9731261196976619662852108486791876458635644729040  \\
\hline
8 & -1.8021253179931622367536330625155209079039086674443  \\
\hline
9 &  3.7133510644596133576858937986178468541115390150895  \\
\hline
10& -11.583138616714443418004214394156033470878899508634  \\
\hline
\vdots & \vdots \\
\hline
20& -7.6931751083769270011123002218244304577221846239268$\times 10^{9}$ \\
\hline
30& -8.1910409909869137068367900925700302382658971132757$\times 10^{20}$ \\
\hline
40& -2.4605043425772457379890548734866774381567481629777$\times 10^{33}$ \\
\hline
50 &-8.9568228254793711194813512752380738598982095960590$\times 10^{46}$ \\
\hline
\end{tabular}
\label{table:nonlin} 
\end{table}

\newpage
\section{Analytical continuation via Mellin transforms}
In the works of Voros [23, p. 683], the application of Mellin transforms yielded two integral representations for $Z(s)$, the first one is

\begin{equation}\label{eq:20}
Z(s)=-\frac{\zeta(s,\tfrac{5}{4})}{2^{s+1}\cos(\frac{\pi}{2}s)}+\frac{1}{\pi}\sin(\frac{\pi}{2}s)\int_{0}^{\infty}\frac{1}{t^s}\left(\frac{\zeta^{\prime}}{\zeta}(\frac{1}{2}+t)+\frac{1}{t-\frac{1}{2}}\right)dt
\end{equation}
which is valid in a narrow strip $0<\Re(s)<1$, and a second transform involving the principal valued integral

\begin{equation}\label{eq:20}
Z(s)=\frac{-\zeta(s,\tfrac{5}{4})+2^{2s}\cos(\pi s)}{2^{s+1}\cos(\frac{\pi}{2}s)}+\frac{1}{\pi}\sin \Big(\frac{\pi}{2}s\Big)\; \text{P.V.}\int_{0}^{\infty}\frac{1}{t^s}\frac{\zeta^{\prime}}{\zeta}(\frac{1}{2}+t)dt
\end{equation}
is actually valid for $\Re(s)<1$. These forms can be obtained from the Chakravarty's functional equation [7]. Both of these integrals are singular at origin, and for (73) there is another pole at $s=\frac{1}{2}$. However, these integrals are still computable that can yield relatively good numerical results. And as for the principal value integral, we consider

\begin{equation}\label{eq:20}
I(s)=\lim_{\epsilon\to 0}\Bigg\{\int_{0}^{\frac{1}{2}-\epsilon}\frac{1}{t^s}\frac{\zeta^{\prime}}{\zeta}(\frac{1}{2}+t)dt+\int_{\frac{1}{2}+\epsilon}^{\infty}\frac{1}{t^s}\frac{\zeta^{\prime}}{\zeta}(\frac{1}{2}+t)dt\Bigg\}
\end{equation}
but the precision in integrating through the zeta function as $\epsilon\to 0$ must be very high, which can be more difficult to compute to high precision. The difficultly gets worse when one tries to compute the derivatives of $Z(s)$ using the principal valued integral, therefore we'll primarily focus on the (73) representation which is valid in $0<\Re(s)<1$, where we will consider computing the derivatives $Z^{(m)}(\frac{1}{2})$ to high precision, so that, based on our previous results (54), we have the full analytical continuation

\begin{equation}\label{eq:20}
Z(s)=\frac{1}{(s-1)^2}\Gamma(\dfrac{s+1}{2})\sum_{n=0}^{\infty}B^{(n)}(\frac{1}{2})\frac{(s-\frac{1}{2})^n}{n!}
\end{equation}
about $a=\frac{1}{2}$, which is valid in $\mathbb{C}\backslash\{1,-1,-3,-5,\ldots\}$. These expansion coefficients are given by

\begin{equation}\label{eq:20}
B^{(n)}(\frac{1}{2})=\lim_{s\to\frac{1}{2}}\Bigg[(s-1)^2\Gamma\Big(\dfrac{s+1}{2}\Big)^{-1}Z(s)\Bigg]^{(n)}
\end{equation}
which from (73) the derivatives are

\begin{equation}\label{eq:20}
Z^{(m)}(\frac{1}{2})=\lim_{s\to\frac{1}{2}}\frac{d^m}{ds^m}\Bigg[-\frac{\zeta(s,\tfrac{5}{4})}{2^{s+1}\cos(\frac{\pi}{2}s)}+\frac{1}{\pi}\sin(\frac{\pi}{2}s)\int_{0}^{\infty}\frac{1}{t^s}\left(\frac{\zeta^{\prime}}{\zeta}(\frac{1}{2}+t)+\frac{1}{t-\frac{1}{2}}\right)dt\Bigg].
\end{equation}
Together, these equations (76), (77), and (78) constitute a full closed-form representation of $Z(s)$ valid in $\mathbb{C}\backslash\{1,-1,-3,-5,\ldots\}$.

We now outline one way that we did to compute these equations numerically. Hence, because there is so many parts to this, a loss of numerical accuracy somewhere is inevitable, especially when computing the derivatives of the integrals. As a result, we'll systematically break it down into small pieces. First, let us split the function into two parts
\begin{equation}\label{eq:20}
Z^{(m)}(s)=\tilde{Z}_1^{(m)}(s)+\tilde{Z}_2^{(m)}(s).
\end{equation}
The first part is the term

\begin{equation}\label{eq:20}
\tilde{Z}_1^{(m)}(s)=\left[-\frac{\zeta(s,\tfrac{5}{4})}{2^{s+1}\cos(\frac{\pi}{2}s)}\right]^{(m)}
\end{equation}
is actually easily numerically computable as well as its derivatives, so we will leave it as it, but second term will be the more complicated one. First we need to separate the integral derivative terms

\begin{equation}\label{eq:20}
\tilde{Z}_2^{(m)}(s)=\left[\frac{1}{\pi}\sin(\frac{\pi}{2}s)I(s)\right]^{(m)}
\end{equation}
by expanding the derivatives of $\tilde{Z}_2^{(m)}(s)$ by the chain rule to generate an expansion

\begin{equation}\label{eq:20}
\begin{aligned}
\tilde{Z}_2^{(m)}(s)=&\frac{(-1)^m}{\pi}\sum_{n=0}^{m}i^{m+n}\left(\frac{\pi}{2}\right)^{m-n}I^{(n)}(s)\binom{m}{n}\times\\
    &\times\Big[i((n+m)\bmod 2)\cos(\frac{\pi}{2}s)+((n+m+1)\bmod 2)\sin(\frac{\pi}{2}s)\Big],
\end{aligned}
\end{equation}
where here, the $(n \bmod 2)$ turns on or off the cosine or sines depending of whether $n$ is even or odd, and the imaginary part gives the correct sign. It generates a workable formula that can be implemented in a program.

Now, the integral to do is
\begin{equation}\label{eq:20}
I^{(m)}(s)=(-1)^m\int_{0}^{\infty}\frac{1}{t^s}\log^m(t)\left[\frac{\zeta^{\prime}}{\zeta}(\frac{1}{2}+t)+\frac{1}{t-\frac{1}{2}}\right]dt,
\end{equation}
and its derivatives, however, there are two major issues that impact computation of higher derivatives of $Z(s)$, one is that at $t=0$ the integrand is improper, and at $t=\frac{1}{2}$ there is a pole cancelation term that has to be integrated through, so with the proceeding analysis, we need to somehow smooth these out.

We first apply the Laurent expansion of the log-zeta function

\begin{equation}\label{eq:20}
-\frac{\zeta^{\prime}(s)}{\zeta(s)}=\frac{1}{s-1}+\sum_{n=0}^{\infty}\eta_n (s-1)^n
\end{equation}
where the eta coefficients can be computed. This series is centered at $s=1$ and has radius of convergence $R=3$ due to the singularity caused by the first trivial zero of $\zeta$. The eta constants are actually very difficult to compute directly, but it's much simpler to express them in terms of the Stieltjes constants by the formula

\begin{equation}\label{eq:20}
\eta_n=(-1)^{n+1}\left[\frac{n+1}{n!}\gamma_n+\sum_{k=0}^{n-1}\frac{(-1)^{k-1}}{(n-k-1)!}\eta_k\gamma_{n-k-1}\right]
\end{equation}
found in Coffey [9, p.532], because the $\gamma_n$ are much more easily computable to high precision and widely available. Now fortunately, the eta series is just enough to be inserted into (83) and cancel the pole, and hence, we obtain the lower integral term

\begin{equation}\label{eq:20}
I_1(t)=\sum_{n=0}^{\infty}\eta_n  \left[\int_{0}^{a}\frac{(-\frac{1}{2}+t)^n}{t^s}\right]dt
\end{equation}
with $1<a<4$ free to choose, but we fix $a=2$ is reasonable value. One could alteratively use an integral identity of the type

\begin{equation}\label{eq:20}
\int_{0}^{a} \frac{(-\frac{1}{2}+t)^n}{t^s}dt  = \frac{2^s(-\frac{1}{2}+t)^{n+1}}{n+1}\: _2F_1\left(n+1,s;n+2;1-2t\right)\bigg\rvert_{0}^{a}
\end{equation}
found by Mathematica [26], which makes a connection to hypergeometric function. From this we define a function for the difference between the anti-derivatives as such

\begin{equation}\label{eq:20}
\delta_n^{(m)}(s)=\frac{\partial^m}{\partial s^m} \left[\frac{2^s}{n+1}\Big((-\frac{1}{2}+a)^{n+1}\: _2F_1\left(n+1,s;n+2;1-2a\right)+(-\frac{1}{2})^{n+1}\: _2F_1\left(n+1,s;n+2;1\right)\Big)\right]
\end{equation}
with its $m^{th}$ derivatives in the $s$ variable in the hypergeometric function is unusually located in the numerator of $_2F_1$, but numerically we find that its  derivatives are smoothly computable, and Mathematica gives these closed-form for these coefficients

\begin{equation}\label{eq:20}
\begin{aligned}
\delta_0^{(0)}(\frac{1}{2})&=2\sqrt{2},\quad \delta_1^{(0)}(\frac{1}{2})=\frac{1}{3}\sqrt{2}, \quad \delta_2^{(0)}(\frac{1}{2})=\frac{23}{15\sqrt{2}} \\
\\
\delta_0^{(1)}(\frac{1}{2})&=\sqrt{2}(4-\log 4), \quad \delta_1^{(1)}(\frac{1}{2})=-\sqrt{2}(\frac{10}{9}+\frac{1}{3}\log 2),
\end{aligned}
\end{equation}
and so on, but we just prefer to compute (87) numerically without resorting to the hypergeometric functions. To sum this up, we get a complete first integral expression as
\begin{equation}\label{eq:20}
I^{(m)}_1(s)=-\sum_{n=0}^{\infty} \eta_{n}\delta_n^{(m)}(s),
\end{equation}
whose main goal is to get rid of the singularity at origin and the pole at $t=\frac{1}{2}$ at high derivatives.

Now, the second integral is already straightforward to integrate and its derivatives are easily found
\begin{equation}\label{eq:20}
I_2^{(m)}(s)=(-1)^{m}\int_{2}^{\infty}\frac{1}{t^s}\log^m(t)\left[\frac{\zeta^{\prime}}{\zeta}(\frac{1}{2}+t)+\frac{1}{t-\frac{1}{2}}\right]dt.
\end{equation}
For higher derivatives, the integrand may seem to be diverging, but eventually it decays to $0$, and so, a note must be taken that numerically the upper limit of integration has to be taken very high. Fortunately, the built-in integrator in Mathematica or Pari can handle it very well. And together, these integrals add up to result in

\begin{equation}\label{eq:20}
I^{(m)}(s)=I^{(m)}_1(s)+I^{(m)}_2(s)
\end{equation}
as required, and valid for $m\geq 0 $. Hence, all in all, this completes the set of equations as an attempt to smoothly compute (78) for high derivatives $m$ without causing any convergence issues.

\begin{table}[hbt!]
\caption{The computation of $B^{(n)}(\tfrac{1}{2})$ (first 50 digits) using the equations (76) to (92). The blue line indicates the last accurate digit} 
\centering 
\begin{tabular}{| c | c |} 
\hline
$\boldsymbol{n}$  & $\boldsymbol{B^{(n)}(\tfrac{1}{2})}$ \textbf{(also in 50 digits)}\\ [0.5ex]
\hline
0 &  0.3160271915014591151297389495924977526309266624541  \\
\hline
1 & -0.3697038073048460995653323747355632929636740870971 \\
\hline
2 &  0.1395921148553775768400923852409921155687771547899  \\
\hline
3 &  0.7936785595326581705505223753582929627453010453660  \\
\hline
4 & -2.7348290310596782885467845906567250844389607748747 \\
\hline
5 &  4.3008854453758201208848420297525288077196862827482 \\
\hline
6 &  2.0066278786744223359450374088654634604561582806775 \\
\hline
7 & -36.962494439621297173325886874777061818393246997054 \\
\hline
8 &  130.48742110411769875204177686724077996724137485496 \\
\hline
9 & -240.92838010377225416324369547429930456498060571889 \\
\hline
10& -87.399353549350560296642724015362132923983838023403 \\
\hline
\vdots & \vdots \\
\hline
20& -157862524.4236541124439435791057859583680249444513 $\times 10^{8}$ \\
\hline
30& -115785913204075.2931497635064477877296431496994334 $\times 10^{14}$ \\
\hline
40& -125665211199394189817.7373673088459255\textcolor{blue}{\underline{7}}444310447888$\times 10^{18}$ \\
\hline
50& -4.78641739986\textcolor{blue}{\underline{0}}6208052121402711032661518765352953191$\times 10^{25}$ \\
\hline
\end{tabular}
\label{table:nonlin} 
\end{table}

We next perform a detailed computation of these equations (76) to (92) by pushing them to the limit, which took over $2$ days, and summarize the results in of $B^{(n)}(\tfrac{1}{2})$ in Table $10$ (also to $50$ digits). To verify these coefficients, we compute the test values

\begin{equation}\label{eq:20}
\begin{aligned}
Z(3)&\approx \frac{1}{4}\sum_{n=0}^{50}B^{n}(\tfrac{1}{2})\frac{(\tfrac{5}{2})^n}{n!}\\
&\approx 0.00072954827270970\underline{4}21587\ldots,
\end{aligned}
\end{equation}
which is accurate to $21$ decimal places. And similarly the value for the $H$ constant

\begin{equation}\label{eq:20}
\begin{aligned}
H &\approx \lim_{s\to 1}\left[ Z(s)-\frac{1}{2\pi(s-1)^2}+\frac{\log(2\pi)}{2\pi(s-1)}-\frac{\log(2\pi)}{4\pi}\right]\\
&\approx -0.01715940430709814\underline{9}96025\ldots,
\end{aligned}
\end{equation}
which is only accurate to $18$ decimal places. We note that using this method, this value was computed independently of the (ADR) [2] or the Brent's method [3][4].  But unlike the (ADR), this series expansion method doesn't assume (RH), but it's much less accurate for about same amount of computational effort. And since we did not characterize the error term, in Table $10$ we note that for $n=40$ and $n=50$ cases, the accuracy is degrading as shown by the last accurate underlined blue digit, and we know this because we compared it with the computation made by the (ADR), which was more accurate.

\section{Conclusion}
The importance of the secondary zeta function $Z(s)$ is that it can analytically represent non-trivial zeros of the Riemann zeta function as $s\to\infty$, such as by equation (13), which is the main motivation behind this article. We develop a new kind of explicit formula relating the secondary zeta function to the prime zeta function, and as a result, we show that the non-trivial zeros can be directly generated from primes in a natural way.

As we developed many series representations for $Z(s)$, we find that the expansion coefficients are difficult to compute numerically, and require special considerations and very high numerical precision, where Pari/Gp software package was instrumental in making these computations. As we have shown several new methods to analytically continue $Z(s)$, but the (ADR) is by far the best one in the terms of achieving high numerical precision. We hope that these results will further aide in the development of many new beautiful formulas and identities in this rather unexplored area, where there should be a lot.

\texttt{Email: art.kawalec@gmail.com}

\section{Appendix A}
The following program implemented in Pari/GP computes the first non-trivial zero by the closed-form formulas (13).

\lstset{language=C,deletekeywords={for,double,return},caption={Pari script for computing first the non-trivial zero using equation (13).},label=DescriptiveLabel,captionpos=b}
\begin{lstlisting}[frame=single]
{
  \\ set limit variable
  m = 20;

  \\ compute the derivative term
  A = 1/(2*m-1)!*derivnum(t=1/2,log(zeta(t)),2*m);

  \\ compute generalized zeta series
  Z = 1/2*(-1)^m*(2^(2*m)-A-1/2^(2*m)*zetahurwitz(2*m,5/4));

  \\ compute the first zero
  t1 = Z^(-1/(2*m));
  print(t1);
}
\end{lstlisting}

\newpage
\section{Appendix B}
The following is a sample code that we developed that implements the (ADR) algorithm [2] in Pari/GP software package [22]. It should generate a few hundreds digits of accuracy for many simple test cases. We have used a fine tweaked version of this code to compute the Tables in this article, but we've used a much bigger and precise database of non-trivial zeros and Bernoulli polynomials, which was obtained by other numerical methods not describe here. Note that precision in Pari must be set very high. This code can also be calibrated to maximize precision, and what we mean by that is to iteratively adjust parameters $a$ and $N$ and comparing $Z(s)$ with known standards near $s$ that we want to compute, such as $Z(2m)$, $Z(0)$ or $Z^{(1)}(0)$, and so on, and this process requires a lot of experimentation. That's how we achieved a significantly higher precision from a calibrated program. 

\lstset{language=C,deletekeywords={for,double,if, return},caption={Pari script for the secondary zeta function using the (ADR) equations (27) thru (32).},label=DescriptiveLabel,captionpos=b}
\begin{lstlisting}[frame=single]
\\ main secondary zeta
Zx(s)=
{
   local(a,A,E,P,S,y);

   \\ set integration limit parameter
   a = 0.005;

   \\ compute the functions
   A=Ax(s,a);
   E=Ex(s,a);
   P=Px(s,a);
   S=Sx(s,a);

   \\ compute Z(s)
   y=A-P+E-S;

   return(y);
}

\\ The A function
Ax(s,a)=
{
  local(tx,y);

  \\ generate a set of non-trivial zeros less than 300
  tx = lfunzeros(1,300);

  y=sum(i=1,40,incgam(s/2,a*tx[i]^2)/gamma(s/2)*1/tx[i]^s);
  return(y);
}

\\ The E function
Ex(s,a)=
{
  1/gamma(s/2)*sum(n=0,300,1/(4^n*factorial(n))*a^(n+s/2)/(n+s/2));
}


\\ The P function
Px(s,a)=
{
  1/(2*sqrt(Pi))*sum(i=2,100,Mangoldt(i)/sqrt(i)
  *incgam((1-s)/2,log(i)^2/(4*a))/gamma(s/2)*(2/log(i))^(1-s));
}

\\ The S function
Sx(s,a)=
{
  local(N,K,A,B,C,y);
  N=100;
  K = a^((s-1)/2)/(4*sqrt(Pi)*gamma(s/2));
  A = -2/(s-1)^2;
  B = (Euler+log(16*Pi^2*a))/(s-1);
  C = 1/gamma(s/2)*1/(8*sqrt(Pi))*sum(i=1,N, 2^(2*i)
  *BernPoly(i,3/4)/factorial(i)*gamma(i/2)*2*a^((s+i-1)/2)/(s+i-1));

  y = K*(A+B)+C;
  return(y);
}

\\ Define von Mangoldt function
Mangoldt(n)=
{
  ispower(n,,&n);
  if(isprime(n),log(n),0)
}

\\ Define Auxiliary function
f(x,t)=t*exp(t*x)/(exp(t)-1);

\\ Define Bernoulli Polynomial
BernPoly(n,x)=
{
  local(delta,y);
  delta = 10^(-500); \\ taking limit
  y = derivnum(t=delta,f(x,t),n);\
  return(y);
}

\end{lstlisting}

\newpage
\section{Appendix C}
For reference, we summarize a high precision computation of the first few odd values of $Z(2n+1)$ to $1000$ decimal digits by the (ADR) algorithm using a highly optimized and calibrated version of the script in Appendix B. The even values can be easily computed by the Voros's formula (3), but there is no known simpler formula for the odd values. We have also computed it to well over $3000$ decimal places, but the new results will be published elsewhere.

\begin{equation}\label{eq:20}
\begin{aligned}
Z(3)=\:& 0.00072954827270970421587551856909397050335150570355423735896527446661230 \\
     & 2447132912878325639671762838465670241435852408263513878826987998494207508 \\
     & 8885085804475186263697199413011480084966272203088728700336243844468256752 \\
     & 3408681057374032516243664750923896178245391471200610831521661864868623771 \\
     & 9981008350821808032694153289752709869393762928727426025438663934785542791 \\
     & 1674351877042641286939396470863595820112331596593164621357724354232059502 \\
     & 5739624542828303360054217223285471146316429365088255140300913431376181876 \\
     & 1120377082063395104378092618263565927716966070075989209338998985346356958 \\
     & 1676943245035728892971601419316025240209233748654242292668468552260654992 \\
     & 9720378193832770070421376800433242194295153056398488540486296586072619189 \\
     & 9163119239368212030984807448938226858830661679960354739515804398696984416 \\
     & 4384126233523215225177568032097452024877692620150229784274171295817361450 \\
     & 7194298180109875520030418160187269226362372525079628507803074554047052409 \\
     & 102077183509495499993543684911252239753259272908325821593956921891770\ldots  \nonumber
\end{aligned}
\end{equation}

\begin{equation}\label{eq:20}
\begin{aligned}
Z(5) =\: &0.00000223118869950210332864062869183719337607643108793448977822617985978 \\
     & 12221524236582470954466136833966440247297286980449542376421449295741935203 \\
     & 69707353398471008236409652933728394361865361430187829399752849921624499279 \\
     & 16060590708347074572906196366778217264313548700998870062035354107184765032 \\
     & 79616935607434826837960593470584989823473238388940832056110254329305367658 \\
     & 52716940584031063832786024218286344325564453145440549670723190640111001097 \\
     & 92959528833041517843917167361065105949171105200217207011958852714776989396 \\
     & 42393573935146268983436444113952793617942300068630586555509789364367446110 \\
     & 17619807347899840742557849020499089533412004004535733523800474071618487076 \\
     & 47250872017729745563528868927469518016098380662133204767575869174392403546 \\
     & 79028332273901424401739847785470262070559006056196112210761072873310156715 \\
     & 53686763553628854778060050195294250651860860233545148729795648690737887072 \\
     & 87097640366826525513747984798038757118898373909779101051927651686645777176 \\
     & 5460646902198601081054720686036774116761436154145645632713172263\ldots \nonumber
\end{aligned}
\end{equation}

\begin{equation}\label{eq:20}
\begin{aligned}
Z(7)=\:&0.000000009675344542702350408719965627461541402301374598418496390368034180494 \\
     & 09539379684057577403986064509491151797773684815935447031528340109149976805166 \\
     & 25356772951431741774963095116941141479776082602903708187984893797113628034854 \\
     & 68268532154205044589073673065905107562612905293340041259438501978075187633929 \\
     & 27471192827610671776970190871522951131875398818469011488991946498613460682248 \\
     & 34404584298387893495130788671879269809486212640632790249021580489963263080369 \\
     & 70552973211207987061707599226050600579678044201891715844948525548491734237877 \\
     & 89131403438039366792658798982631395800414034492124755789490732393268322848249 \\
     & 58759893932928887685381405053476092418802847575434128157349311125790029085018 \\
     & 29560900708672613699416652806107469391901918922572783448775000383665026223266 \\
     & 26014270542205662930014290872841168715803750103657252655060670167848893578301 \\
     & 10633804334384463550107213669830002448886753769212496201843916386837590112504 \\
     & 73881659154345592066274061514667037586810080799129631434987796858520634167166 \\
     & 429726297412562471999519416\ldots \nonumber
\end{aligned}
\end{equation}

\begin{equation}\label{eq:20}
\begin{aligned}
Z(9)=\:&0.000000000045991912392894862969274651949689848934500247171975676359515259476 \\
      & 17226302282722719145458884398086960343237235698529918246425420830851460628161 \\
      & 44818931898092840111991568905627494946673555987116556295643670981922579742233 \\
      & 25017792625552834429279255898637621341474655917544411907134138321244373074475 \\
      & 02355482364221240268070144432982679112045225478991739794118959787677747399405 \\
      & 76062488421634707537437210624515786098603345347151197406381555055885359462075 \\
      & 24155957863639566607642218075762559492030752657158634598285202555417909859283 \\
      & 22263862265823678369012741616541691632771436544856239973444318215950791974138 \\
      & 74891796897732450600803290705541548951813520594420260225873710779351860802674 \\
      & 82684509984937793968364427781166697886532091703224544425294182020653828733282 \\
      & 19684764553150426711864351049756617936616147018172181663853725934524153054424 \\
      & 43707269217732260973945250219212833930647796033494112797829468559193886546273 \\
      & 53523878848316078130956364043307862245352111573749108563814339444492039484561 \\
      & 43640752017433973681267621121\ldots \nonumber
\end{aligned}
\end{equation}

\begin{equation}\label{eq:20}
\begin{aligned}
Z(11)=\:&0.000000000000225565062515596642898367454188027063536645929510183777475351311 \\
       & 11543488351467458681563430243716773307332687818692894694459742051046453109827 \\
       & 10931603716573630012805945670157902717927060524850998719989294914141843515475 \\
       & 64392766591816696308533079294240809779029553029395441199345217372823641107765 \\
       & 67214306571979662283433359155691470141745266085834256781066510518318669115511 \\
       & 96666211341908201367200023233967123252901254098862416560881040838067888142021 \\
       & 70952792990355223063662015657847860282511831545777699249088652431401193766899 \\
       & 92671292737441996399149758401885912751199763513005354667955189070388041573689 \\
       & 52248022167305733732890059873790706088680801374144727857674923675190990611145 \\
       & 73683045883403234812426540644202639320742907740070982079881985178315433565646 \\
       & 86648126419298040509311902780712795190301456111955904338501262751562909100149 \\
       & 66723306435572281960378588791009666128964082659222463750654601452537092482300 \\
       & 36502551457296218087907850537504726201334954904176395019870856887496953459139 \\
       & 505495805823573937334659950273\ldots \nonumber
\end{aligned}
\end{equation}

\end{document}